\documentclass[11pt,a4paper]{amsart}

\usepackage[latin1]{inputenc} \usepackage[T1]{fontenc}

\usepackage{graphicx}

\usepackage{multicol, psfrag, subfigure}

\usepackage{amssymb}
\usepackage{amsthm}

\usepackage[dvipsnames]{xcolor}
\xdefinecolor{forestgreen}{named}{ForestGreen}

\usepackage[all]{xy}
\usepackage{stmaryrd}

\everymath{\displaystyle}

\newtheorem{thm}{Theorem}[section] 
\newtheorem{prop}[thm]{Proposition}
\newtheorem{df}[thm]{Definition} 
\newtheorem{lem}[thm]{Lemma}
\newtheorem{cor}[thm]{Corollary}
\newtheorem{ex}[thm]{Example}
\newtheorem{rk}[thm]{Remark}

\newcommand{\R}{\mathbb{R}}
\newcommand{\Q}{\mathbb{Q}}
\newcommand{\dd}{\mathrm{d}}
\newcommand{\T}{\mathbf{T}}
\newcommand{\p}{\underline{+}}
\newcommand{\ti}{\underline{\times}}
\newcommand{\s}{\underline{\square}}

\newcommand{\Pv}{\noindent \textit{Proof.}}

\begin{document}             

\date{May 2017}

\bibliographystyle{plain}

\title[Operations on Legendrian submanifolds]{Operations on Legendrian submanifolds}

\newcommand{\cs}{$^\dagger$} \newcommand{\cm}{$^\ddagger$}

\author{Maÿlis Limouzineau}

\newcommand{\nfont}{\fontshape{n}\selectfont}

\address{Institut de mathématiques de Jussieu-Paris-rive-gauche\\ Université Pierre et Marie Curie (Paris 6).}

\begin{abstract}
We focus on Legendrian submanifolds of the space of one-jets of functions, $J^1(\R^n,\R)$. We are interested in processes - operations - that build new Legendrian submanifolds from old ones. We introduce in particular two operations, namely the sum and the convolution, which in some sense lift to $J^1(\R^n,\R)$ the operations sum and infimal-convolution on functions that belong to convex analysis. We show that these operations fit well with the classical theory of generating functions. Finally, we refine this theory so that the min-max selector of generating functions plays its natural role.  
\end{abstract}

\maketitle

\section*{Introduction}

\vspace{0.2cm}

\indent The present paper concerns Legendrian submanifolds of spaces of $1$-jets of functions equipped with their standard contact structure. All our objects -- submanifolds, applications or functions -- are assumed to be smooth if nothing else is specified. Let $n$ be an integer. We study three operations on Legendrian submanifolds of $J^1(\R^n,\R)$: the \textit{sum}, the \textit{convolution} and the \textit{transformation} $\mathbf{T}$. By "operations" we mean procedures which build new Legendrian submanifolds from old ones.

The $1$-jet space $J^1(\R^n,\R)=\R \times T^* \R^n$ is endowed with coordinates $(u,q,p)$, where $q$ stands for the base space of $T^*\R^n$ and $p$ corresponds to fibers $T^*_q\R^n$. The standard contact structure of $J^1(\R^n,\R)$ is given as the kernel of the differential $1$-form $\dd u - p\dd q$. 

Let $L_1$ and $L_2$ be two Legendrian submanifolds of $J^1(\R^n,\R)$. The sum $L_1 \p L_2$ is defined by making the sum in each $(u,p)$-fiber, $$L_1 \p L_2 = \lbrace (u_1+u_2,q,p_1+p_2) \vert (u_i,q,p_i)\in L_i \rbrace.$$ 

On the other hand, as the base of the vector bundle $T^*\R^n$ is vectorial, it allows to make the sum in the base space rather than in the $p$-fibers. The convolution $L_1 \s L_2$ is then defined by $$L_1 \s L_2 = \lbrace (u_1+u_2,q_1+q_2,p) \vert (u_i,q_i,p)\in L_i \rbrace.$$ 

Those two operations are naturally linked by a Legendre transform of $J^1(\R^n,\R)$, which exchanges the roles of the base and the $p$-fibers. In this note we work with the transformation 
\begin{align*}
\T~: \ & J^1(\R^n,\R) \rightarrow J^1(\R^n,\R)\\
& (u,q,p) \mapsto (qp-u,p,q).
\end{align*}

\begin{thm}[Theorem \ref{thm21}]
If $L_1$ and $L_2$ are two Legendrian submanifolds of $J^1(\R^n,\R)$, their sum $ L_1 \p L_2 $ and their convolution $L_1 \s L_2$ are linked by the following identities: $$\underline{i}) \ \mathbf{T}( L_1 \p L_2 )=(\mathbf{T}L_1) \s (\mathbf{T}L_2) \quad \text{ and } \quad  \underline{ii}) \  \mathbf{T}(L_1 \s L_2 )=(\mathbf{T}L_1) \p (\mathbf{T}L_2 ).$$ 
\end{thm}

Beyond introducing these three operations on Legendrian submanifolds, our purpose is to show how they correspond to the lifts of three classical operations on functions coming from the domain of convex analysis. Remind that the first examples of Legendrian submanifolds of $J^1(\R^n,\R)$ are given by functions: consider $f$ a function defined on $\R^n$, its \textit{$1$-graph}  
$$j^1 f = \lbrace (u,q,p) \ \vert \ q \in \R^n , u=f(q), p=\frac{\partial f}{\partial q}(q)\rbrace $$
is a Legendrian submanifold of $J^1(\R^n,\R)$. We view the \textit{wave fronts} of Legendrian submanifolds -- i.e. their projection on the space of $(u,q)$-coordinates -- as graphs of generalized (multivalued) convex functions. Then the sum, the transformation $\T$ and the convolution are the natural generalizations of the lifts of three operations on (single-valued) functions: the sum, the \textit{Legendre--Fenchel transform} and the \textit{infimal-convolution} (see also \cite{1}).\\ 

The idea of a reconciliation between these operations on functions used in convex analysis and contact phenomena is not new. The Legendre--Fenchel transform and the infimal-convolution are commonly used in the field of thermodynamics (see for instance \cite{17}), and the links between thermodynamics and contact setting are emphasized a number of times in the mathematics and physics literature (for instance recently in \cite{16}). Without going for details, just remind the fundamental thermodynamic relation $$\dd U=T \dd S -P \dd V \ ,$$
 where $U$ stands for the internal energy, $T$ the temperature, $S$ the entropy, $P$ the pressure and $V$ the volume. Roughly speaking, it follows that the set of stable states for a closed thermodynamical system is a Legendrian submanifold of $\R^5$.\\ 

If $f_1$ and $f_2$ are two functions defined on $\R^n$, the sum $f_1 + f_2$ is also a function defined on $\R^n$, which inherits of the regularity of $f_1$ and $f_2$ (if $f_1$ and $f_2$ are $C^{\infty}$ for instance, so does $f_1 +f_2$). The sum operation applied on the $1$-graphs of $f_1$ and $f_2$ clearly correspond to the $1$-graph of the sum $f_1+f_2$. In that sense, our sum operation lifts the sum of functions to $J^1(\R^n,\R)$, and extends it to Legendrian submanifolds.

Let $f$ be a function defined on $\R^n$ such that its gradient $\nabla f$ defines a diffeomorphism of $\R^n$. The \textit{Legendre transform} associates to $f$ another smooth function defined on $\R^n$ whose gradient is also a diffeomorphism, equal to $(\nabla f)^{-1}$. The Legendre--Fenchel transform permits to extend the Legendre transform to functions which do not have an invertible gradient, but are \textit{almost-convex} (Definition \ref{def37}). The Legendre--Fenchel transform of $f$ is the function $$f^*: \R^n \rightarrow \R \cup \lbrace +\infty \rbrace \ , \ q \mapsto \sup_{v \in \R^n}\lbrace qv -f(v)\rbrace ,$$
where $qv$ stands for the standard Euclidean scalar product. When $f$ is convex with its gradient $\nabla f$ invertible, so is $f^*$. In that case the $1$-graph of $f^*$ coincides with the transformation $\T$ of the $1$-graph of $f$. The Legendre--Fenchel transform has no interest if the function $f$ does not have an affine minimizer -- $f$ is not \textit{admissible} in the sense of \cite{17}. If a function is almost-convex, the Legendre--Fenchel transform gives back a convex -- but not necessarily smooth -- function. Applying the Legendre--Fenchel transform twice gives the convex approximation of $f$. 

Let $f_1$ and $f_2$ be two functions on $\R^n$, the infimal-convolution of $f_1$ and $f_2$ is the function
$$f_1 \square f_2 : \R^n \rightarrow \R \cup \lbrace -\infty \rbrace \ , \ q \in \R^n \mapsto \inf_{v \in \R^n}\lbrace f_1(v)+f_2(q-v)\rbrace \ . $$

Like the Legendre--Fenchel transform, the infimal-convolution is relevant only for almost-convex functions. If $f_1$ and $f_2$ are convex and have invertible gradients, so is $f_1 \square f_2$ and its $1$-graph coincides with the convolution of the $1$-graphs of $f_1$ and $f_2$.  \\
 
Remarkably, the trace of the Fourier-type identities which come naturally at the geometrical level of Legendrian submanifolds (Theorem \ref{thm21}) appears here. Indeed it is part of the classic convex analysis that for convex or admissible functions \cite{1}, \cite{17} -- we have $$i) \ (f_1 \square f_2)^* =f_1^*+f_2^* \quad \text{and} \quad ii) \ (f_1+f_2)^*=f_1^* \square f_2^* \ .$$

The geometrical operations of transformation $\T$ and convolution do not exactly correspond to the lifts and generalisations of these operations on functions at the level of Legendrian submanifolds. Looking at almost-but-not-convex functions, the transformation $\T$ and the convolution of the associated $1$-graphs do not produce $1$-graphs of functions. The wave fronts of the results are not single-valued graphs but multi-valued ones. From the wave front of the transformation $\T$ of $j^1f$, the definition by supremum in the Legendre--Fenchel transform of $f$ continuously selects a single value above each point $q$ of the base space. Similarly, the definition by infimum in the infimal-convolution of $f_1$ and $f_2$ gives a continuous section of the wave front of the convolution of $j^1f_1$ and $j^1f_2$.
To refine our point of view and the comparison between the three operations on functions (sum, Legendre--Fenchel transform and infimal-convolution) and the three operations on Legendrian submanifolds (sum, transformation $\T$ and convolution), we have to consider the \textit{generating function} setting to use an important object in symplectic geometry: the \textit{min-max selector}. 

The min-max selector was originally introduced by Sikorav and Chaperon to define weak solutions of the Hamilton--Jacobi equation \cite{MR1094198}, \cite{MR2276957}. Its vocation is to single out a continuous sub-graph in the wave front of Legendrian submanifolds realized by some generating functions. We view convex functions as \textit{simple functions} (Definition \ref{def23}) with \textit{Morse index} (Definition \ref{def35}) equal to zero. The selector allows to extend the selection procedure for the class of \textit{almost-simple} functions defined on $\R^n$ of any Morse index in $\llbracket 0,n \rrbracket$.

In the \textit{almost-convex} (resp. \textit{almost-concave})\footnote{See Definition \ref{def37}.} case, it is known that the min-max selects the minimum (resp. the maximum). We recover the Legendre--Fenchel transform and the infimal-convolution for almost-convex functions using the selector of generating functions (Corollary \ref{cor31}). 

Legendrian submanifolds realized by generating functions constitute an intermediate class of Legendrian submanifolds, between the set of $1$-graphs of functions and the set of all possible Legendrian submanifolds of $J^1(\R^n,\R)$. The three operations on Legendrian submanifolds behave remarkably well with respect to the notion of generating function (Lemma \ref{lem31}). How do the three operations on Legendrian submanifolds interact in general with the selector of generating functions? One of the objectives is to recover the trace of the Fourier-type identities at the level of the selector. We get from equivalent generating function theory the analogous identities of $i)$ and $ii)$ for the selector (Theorem \ref{thm327}). They generalize $i)$ and $ii)$, and confirms that these identities are not a specificity of the convex setting, but a generating function manifestation.\\   

The first section is devoted to an elementary Legendrian geometry tool box. Most of it is folklore. We take the opportunity to remind the definitions of three classical operations on Legendrian submanifolds: the \textit{product}, the \textit{slice} and the \textit{contour}. They can be used to recover the operations sum and convolution of Legendrian submanifolds (Remark \ref{rk21}). Section $2$ is devoted to the operations sum and convolution and their properties. Finally in section $3$ we work on generating functions and the notion of min-max selector. We remind the classical notion of generating function for Legendrian submanifolds, and associate to each operation on Legendrian submanifolds the corresponding operation on generating functions (Lemma \ref{lem31}). Then, we search for each operation a class of generating functions that admits a selector and is stable, so that the selector still exists at the end of the operation. We introduce adapted classes of generating functions for operations sum and transformation $\T$, named \textit{almost-simple} (Definition \ref{def39}) -- which is slightly more general than the usual \textit{quadratic at infinity} notion \cite{4} -- and \textit{globally almost-simple} (Definition \ref{def312}).

Unfortunately, this last class of generating functions is not stable by sum, unless the Morse index is set to be minimal or maximal. Restricting to the class of almost-convex functions (when the Morse index is minimal equal to zero) allows to work with the three operations together with a persisting selector. There we recover the convex analysis setting. 

A maximal class of generating functions stable for the three operations at the same time is still to be specified.
 
\section*{Acknowledgements}
I am grateful to my PhD advisor, Emmanuel Ferrand, for his guidance all along this work. I also want to thank the members of ANR COSPIN, in particular Sheila Sandon for numerous pieces of advice, as well as François Laudenbach, and finally Valentine Roos for enlightening conversations concerning the min-max selector.   

\vspace*{1cm}
\section{Preliminaries: Legendrian things in $J^1(\R^n,\R)$}\label{sec1}

\vspace{0.2cm}

Throughout this paper, we make the assumption that our functions take finite values, unless something else is specified.

\subsection{Contact structure} We work in the space $J^1(\R^n,\R)=\R \times T^*\R^n$ of $1$-jets of functions defined on $\R^n$, for $n \in \mathbb{N}^*$. Moreover, we use the standard euclidian scalar product to systematically identify $T^*\R^n$ with the product $\R^n \times \R^n$ and use canonical coordinates $(q_1,\ldots,q_n,p_1, \dots ,p_n)$. So at the end we work with coordinates $(u,q,p)=(u,q_1,\ldots,q_n,p_1,\ldots,p_n)$ on $J^1(\R^n,\R)$, viewed as the product space $\R \times \R^n \times \R^n$. We will refer to the space of $q$'s as \textit{the base space}.

\begin{df}\label{def11} Consider the differential $1$-form $\mathrm{d}u-p\mathrm{d}q=\mathrm{d}u - \sum_{i=1}^np_i \mathrm{d}q_i$ on $J^1(\R^n,\R)$.
The \textbf{standard contact structure} on $J^1(\R^n,\R)$ is the hyperplane field $\xi$ defined as the kernel of this $1$-form: $$\xi=\text{ker}(\mathrm{d}u-p\mathrm{d}q) \ . $$
\end{df}

\noindent \textbf{Notations} \\ $\bullet$ \ If $f$ is a function defined on $\R^n$, $\nabla f$ denotes the gradient of $f$: \begin{align*}
\nabla f \quad : \quad \R^n \quad & \longrightarrow \quad \R^n \\
q \quad & \longmapsto \left( \frac{\partial f}{\partial q_j}(q)\right)_{j=1,\ldots , n} 
\end{align*}
$\bullet$ \ For $q$ and $v$ in $\R^n$, $qv$ denotes their Euclidean scalar product: \ $qv=\sum_{i=1}^n q_iv_i$, \ $q=\left( q_1,\dots , q_n\right)$ , $v=\left(v_1, \dots , v_n \right)$.\\ 

\begin{df}\label{def12}
A \textbf{Legendrian submanifold} $L \subset J^1(\R^n,\R)$ is a $n$-dimensional submanifold, which is everywhere tangent to the contact structure. In other words: $$(\mathrm{d}u-p\mathrm{d}q)_{\vert L}\equiv0 \ . $$ 
\end{df}

\begin{ex}\label{ex11}
Let $f$ be a function defined on $\R^n$, its \textbf{$1$-graph} is $$j^1f=\lbrace ( u=f(q),q,p=\nabla f (q) ) \ \vert \ q \in \R^n \rbrace \ ,$$
This is the simplest example of Legendrian submanifolds of spaces of $1$-jets of functions. Conversely, if a Legendrian submanifold is a graph over the base $\R^n$ in $J^1(\R^n,\R)$, one can show that it must take the form of a $1$-graph of function.\\
\end{ex}

\begin{rk} The cotangent bundle $T^*\R^n$ is naturally endowed with the standard symplectic form $\omega=\dd( p \dd q)$.
The projection of a Legendrian submanifold $L \subset J^1(\R^n,\R)$ on $T^*\R^n$ is an immersed exact Lagrangian submanifold $\mathcal{L}$. 
\end{rk}

If a Legendrian submanifold $L$ is in generic position, all its description can be recovered from $(u,q)$'s coordinates by defining the missing coordinates $p=(p_1,\ldots,p_n)$ as the slopes $\left(\frac{\partial u}{\partial q_1},\ldots, \frac{\partial u}{\partial q_n}\right)$. It permits to work on Legendrian submanifolds using the projection on the space of $0$-jets of functions\footnote{This allows us to draw pictures for Legendrian submanifolds of dimensions $1$ and $2$.}: \begin{align*}
\mathbf{pr} \ : \ J^1(\R^n,\R)=\R \times \R^n \times \R^n &\longrightarrow \ J^0(\R^n,\R)=\R \times \R^n\\
(u,q,p) \quad & \longmapsto \quad \quad (u,q) 
\end{align*}

\begin{df}\label{def13} 
For $L$ a Legendrian submanifold, $\mathbf{pr}(L)$ is called the \textbf{wave front of $L$}.\\ 
\end{df}

The $1$-graph of a function $f$ defined on $\R^n$ simply projects on the graph of $f$: $$\mathbf{pr}(j^1f)=\lbrace (f(q),q) \ \vert \ q \in \R^n \rbrace \ . $$ 
More generally, generic Legendrian submanifold $L \subset J^1(\R^n,\R)$ is almost everywhere locally a graph over the base space. It projects on an $n$-dimensional object of $J^0(\R^n,\R)$, which is almost everywhere a $n$-dimensional submanifold, away from a singular subset of codimension greater than $0$. 

\begin{ex}\label{ex12}
\ $\bullet$ \ In the case $n=1$, two types of singularity may appear generically in a wave front: double points and (right or left) cusps.\\
$\bullet$ \ In the case $n=2$, a "swallow tail" may appear, as drawn if Figure \ref{fig2}, as well as lines of double points, lines of cusp, \dots (see \cite{2} for the exhaustive list of generic two dimensional local wave fronts).\\
\end{ex}

\begin{df}\label{def14}
A \textbf{contactomorphism} of $J^1(\R^n,\R)$ is a diffeomorphism of $J^1(\R^n,\R)$, which preserves the contact structure $\xi$, i.e. such that the standard contact form $\mathrm{d}u-p\mathrm{d}q$ is sent to itself modulo multiplication by a smooth nowhere vanishing function.\\
\end{df}

\begin{ex}\label{ex13}
The map \begin{align*}
\mathbf{T} :\quad  & J^1(\R^n,\R) \longrightarrow \ J^1(\R^n,\R)\\
& \quad (u,q,p) \longmapsto (pq-u,p,q)
\end{align*}
is a contactomorphism. We refer to it as the \textbf{transformation $\T$}.\\ 
\end{ex}

\begin{figure}[h]
\begin{center}
\includegraphics[scale=0.65]{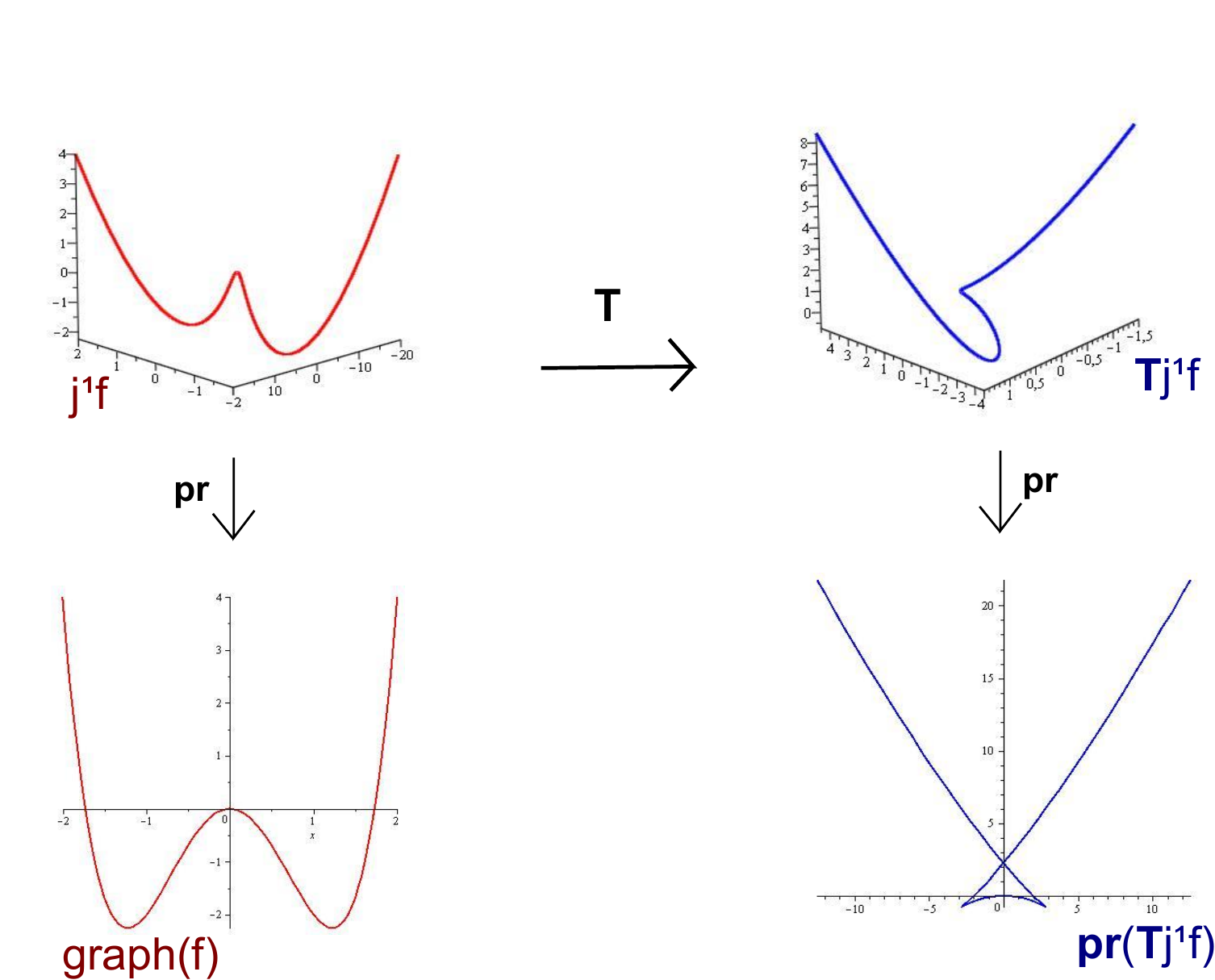}
\end{center}
\caption{$\T$ applied to the $1$-graph of $f : q \mapsto q^4-3q^2$}
\label{fig1}
\end{figure}

As it preserves the contact structure, a contactomorphism transforms Legendrian submanifolds into Legendrian submanifolds. If one starts from the $1$-graph of some function $f$, that is $p=\nabla f (q)$, then $\mathbf{T}j^1 f$ is another Legendrian, but it is not necessarily the $1$-graph of a function any more. If not, its front $\textbf{pr}\left(\mathbf{T}j^1 f \right)$ is not the graph of a function, but a multivalued graph -- see figure \ref{fig1} which describes $\mathbf{T}$ applied to the $1$-graph of $f : \R \rightarrow \R \ , \ q \mapsto q^4-3q^2$, viewed in $J^1(\R^n,\R)$ and in $J^0(\R^n,\R)$ after the front projection.

\subsection{Product} Let $n_1$ and $n_2$ be natural integers. From two Legendrian submanifolds $L_1 \subset J^1(\R^{n_1},\R)$ and $JL_2 \subset J^1(\R^{n_2},\R)$, we naturally define a Legendrian submanifold of $J^1(\R^{n_1} \times \R^{n_2},\R)$. Note that it is not the product, in the set-theoretic sense, of $L_1$ by $L_2$, which sit in $J^1(\R^{n_1},\R)\times J^1(\R^{n_2},\R)$.
\begin{df}\label{def15}
Consider the spaces $J^1(\R^{n_1},\R)$ and $J^1(\R^{n_2},\R)$ with respective coordinates $(u_1,q_1,p_1)$ and $(u_2,q_2,p_2)$, and then the coordinates $(u,q_1,q_2,p_1,p_2)$ in $J^1(\R^{n_1+n_2},\R)$. Given $L_1 \subset J^1(\R^{n_1},\R)$ and $L_2 \subset J^1(\R^{n_2},\R)$, the \textbf{product} $L_1 \ti L_2$ is defined by
\vspace{0.2cm}
\begin{center}
$L_1 \ti L_2 =\lbrace (u_1+u_2,q_1,q_2,p_1,p_2) \ \vert \ (u_i,q_i,p_i) \in L_i \rbrace$.\end{center}
\end{df}

\vspace{0.1cm}
\begin{lem}\label{lem111} If $L_1$ and $L_2$ are Legendrian submanifolds, then $L_1 \ti L_2$ is an immersed Legendrian submanifold.
\end{lem}

\Pv \ Consider $\mathcal{L}_1 \subset T^* \R^{n_1}$ and $\mathcal{L}_2 \subset T^* \R^{n_2}$ the corresponding Lagrangian projections of $L_1$ and $L_2$. They are immersed Lagrangian submanifolds. $L_1 \ti L_2$ is a lift of the (usual) product $\mathcal{L}_1\times \mathcal{L}_2 \subset T^* (\R^{n_1} \times \R^{n_2})$. Thus $L_1 \ti L_2$ is an immersed submanifold of dimension $n_1+n_2$. The standard contact form $\mathrm{d}u-(p_1\mathrm{d}q_1+p_2\mathrm{d}q_2)$ vanishes on it, so $L_1 \ti L_2$ is Legendrian.\qed \\

In \cite{20}, P. Lambert-Cole studies a more general case of Legendrian product. He notices that if the sets of Reeb actions (i.e. lenghts of Reeb chords) of $L_1$ and $L_2$ are disjoint, then the product is embedded. It shows in particular that $L_1 \ti L_2$ is generically embedded. 
%\textcolor{red}{*** c'est peut être maladroit dit comme ça ***}

\begin{ex}
We depict in Figure \ref{fig03} an example of a product of two trivial Legendrian knots. The result is the embedded Legendrian surface of $J^1(\R^2,\R)$, whose wave front in depicted in Figure \ref{fig03}. Let us call this wave front a \textit{toric pillow}, with four corners which are non-generic singularity of wave fronts (see the classification in \cite{2}).

\begin{figure}[ht]
\begin{center}
\includegraphics[scale=0.5]{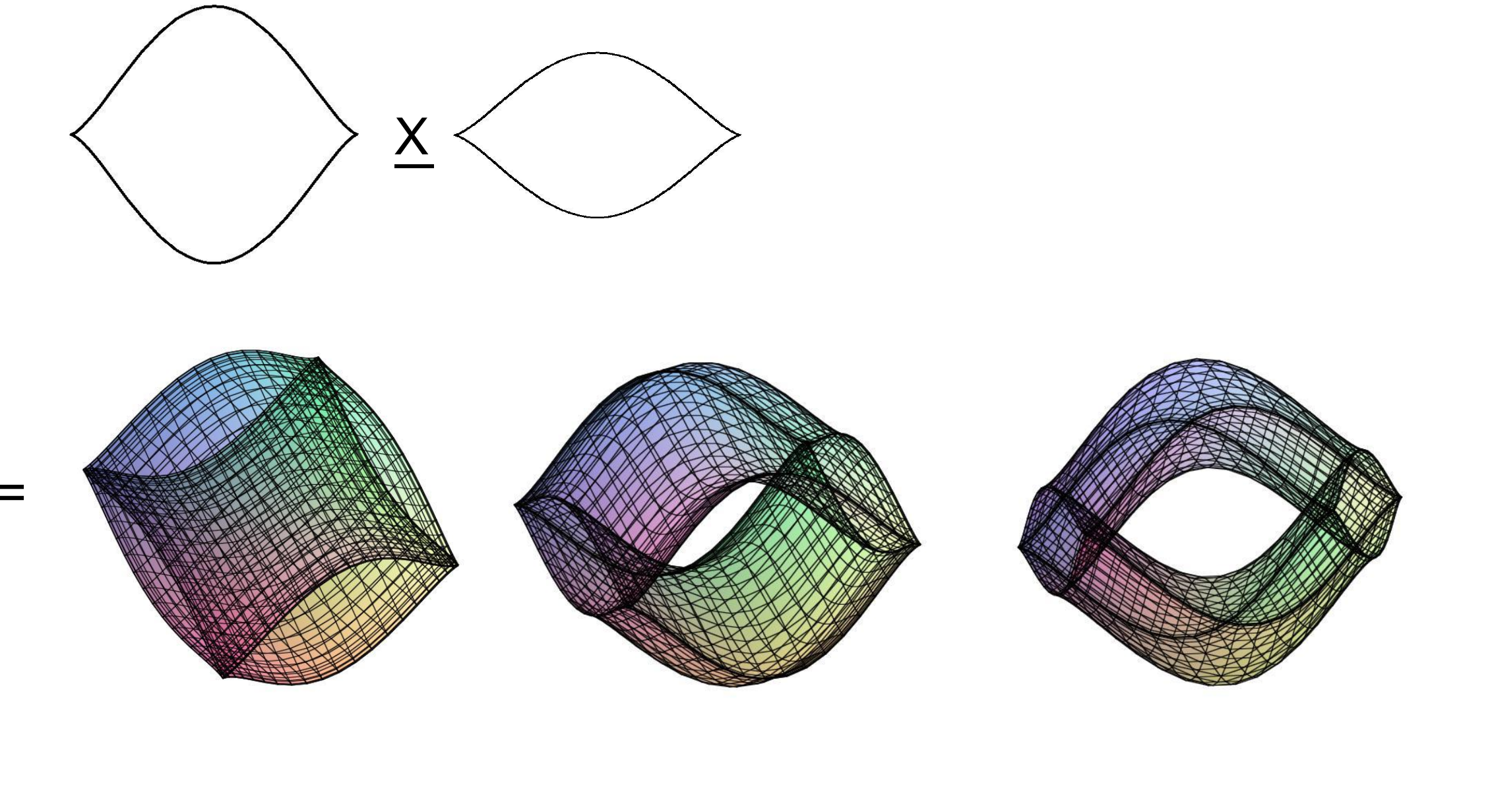}
\end{center}
\caption{A \textit{toric pillow} (the result is presented with three different angles of view) presents four corners, each corresponding to a $D_4^*$ singularity.}
\label{fig03}
\end{figure}
\end{ex}

\begin{rk}\label{rk12}
A direct computation shows that the transformation $\mathbf{T}$ commutes with the product: $$\mathbf{T}(L_1\ti L_2) =(\mathbf{T}L_1)\ti (\mathbf{T}L_2) \ .$$

\noindent Note that $\mathbf{T}$ is applied here in three different spaces: $J^1(\R^{n_1+n_2},\R)$, $J^1(\R^{n_1},\R)$, or $J^1(\R^{n_2},\R)$. However we will keep the same notation $\mathbf{T}$ whichever space is concerned.\\
\end{rk}

\subsection{Slice and Contour} These are also called operations of reduction, as they permits to built "small" Legendrian submanifolds of $J^1(\R^{n_1},\R)$ by restricting "big ones" of $J^1(\R^{n_1} \times \R^{n_2},\R)$. 

To define the first one, we consider the following subspace of $J^1(\R^{n_1+n_2},\R)$: $$ \mathcal{E}_{\sigma}=\lbrace (u,q_1,q_2,p_1,p_2) \ \vert \ q_2=0 \rbrace \ .$$ There is a natural projection from $\mathcal{E}_{\sigma}$ to $J^1(\R^{n_1},\R)$, which forgets the dual coordinate $p_2$, that is
\begin{align*}
\mathbf{p}_{\sigma}\quad : \qquad \qquad \mathcal{E}_{\sigma} \quad \quad \quad & \longrightarrow \ J^1(\R^{n_1},\R)\\
 (u,q_1,0,p_1,p_2) & \longmapsto  \ (u,q_1,p_1)\\
\end{align*}

\begin{df}\label{def16}
Let $L$ be a subset of $J^1(\R^{n_1+n_2},\R)$. The \textbf{slice of $L$ along $\R^{n_1} \times \lbrace 0 \rbrace$} is the set $\mathbf{p}_{\sigma}(L \cap \mathcal{E}_{\mathbf{\sigma}})\subset J^1(\R^{n_1},\R)$, denoted $\sigma(L)$. 
\end{df}

\begin{figure}[ht]
\begin{center}
\includegraphics[scale=0.5]{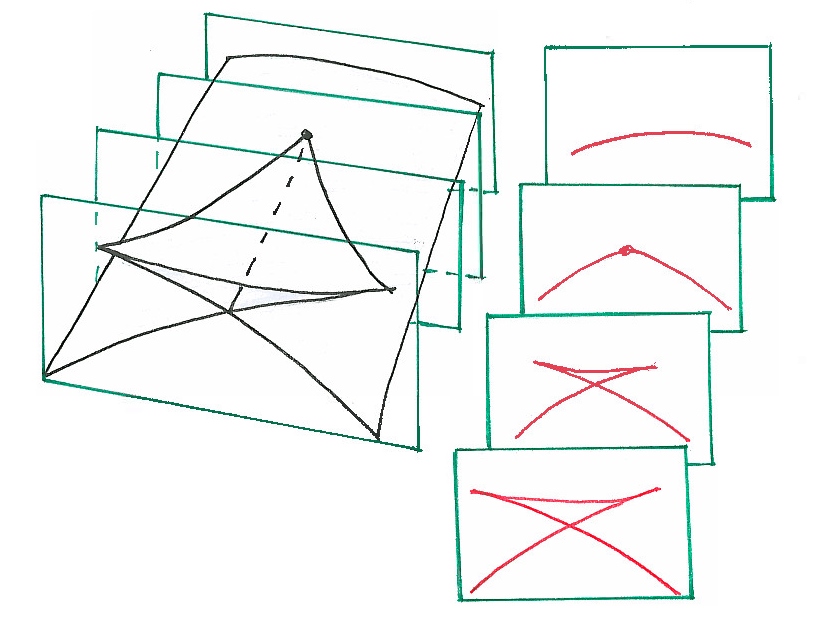}
\end{center}
\caption{Some slices of the swallow tail at the level of fronts}
\label{fig2}
\end{figure} 

In the same way, a dual operation is built by considering the subspace of $J^1(\R^{n_1+n_2},\R)$ $$ \mathcal{E}_{\kappa}=\lbrace (u,q_1,q_2,p_1,p_2) \ \vert \ p_2=0 \rbrace \ .$$ We denote by $\mathbf{p}_{\kappa}$ the projection which consists in forgetting the $q_2$ coordinate
\begin{align*}
 \mathbf{p}_{\kappa}\quad : \qquad \qquad \mathcal{E}_{\kappa} \quad \quad \quad & \longrightarrow J^1(\R^{n_1},\R)\\
 (u,q_1,q_2,p_1,0) & \longmapsto (u,q_1,p_1)
\end{align*} 

\begin{df}\label{def17}
Let $L$ be a subset of $J^1(\R^{n_1+n_2},\R)$. The \textbf{contour of $L$ in the direction of $\R^{n_1} \times \lbrace 0 \rbrace$} is the set $\mathbf{p}_{\kappa}(L \cap \mathcal{E}_{\kappa})\subset J^1(\R^{n_1},\R)$, noted $\kappa(L)$. 
\end{df}

\begin{figure}[ht]
\begin{center}
\includegraphics[scale=0.6]{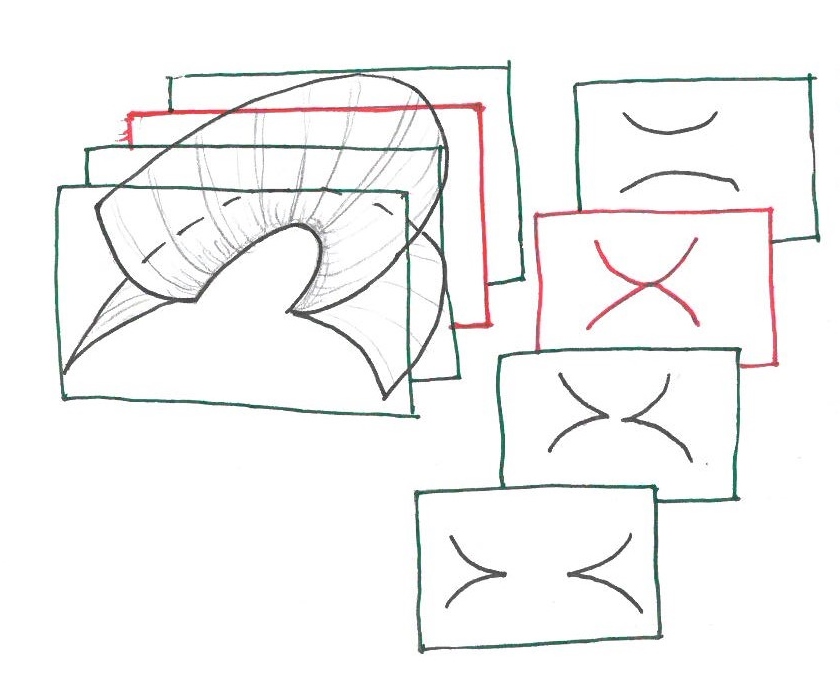}
\end{center}
\caption{Slices of a curved line of cusps gives rise to a degenerate case. The third one correspond to a non-immersed Legendrian}
\label{fig3}
\end{figure}

\begin{lem}\label{lem11} \ Let $L$ be a Legendrian submanifold of $J^1(\R^{n_1+n_2},\R)$. If $L$ is transverse to $\mathcal{E}_{\sigma}$ (respectively $L$ is transverse to $\mathcal{E}_{\kappa}$), then the slice $\sigma(L)$  (respectively the contour $\kappa(L)$) is an immersed Legendrian submanifold of $J^1(\R^{n_1},\R)$.
\end{lem}

\Pv \ The projection of the Legendrian submanifold $L \subset \R \times T^*\R^{n_1 +n_2}$ on $T^*\R^{n_1+n_2}$ is an exact Lagrangian immersion in $T^*\R^{n_1+n_2}$ endowed with the standard symplectic form: \ $\dd p_1 \wedge \dd q_1 + \dd p_2 \wedge \dd q_2 \ .$ For all $x=(u,q_1,q_2,p_1,p_2) \in L$, the coordinates of two tangent vectors $X$ and $Y$ in $T_x L$ must satisfy the equation $$(\star) \quad \quad X_{p_1}Y_{q_1}-Y_{p_1}X_{q_1}+X_{p_2}Y_{q_2}-Y_{p_2}X_{q_2}=0  \ ,$$
with $X=(X_u,X_{q_1},X_{q_2},X_{p_1},X_{p_2})$ and $Y=(Y_u,Y_{q_1},Y_{q_2},Y_{p_1},Y_{p_2})$ in  $$T_x(J^1(\R^{n_1+n_2},\R)) \simeq \R \times \R^{n_1} \times \R^{n_2} \times \R^{n_1} \times \R^{n_2} \ . $$

Let us assume that $L \pitchfork \mathcal{E}_{\sigma}$. Then, $L \cap \mathcal{E}_{\sigma}$ is a submanifold in $J^1(\R^{n_1+n_2};\R)$, and
$$\forall x \in L \cap \mathcal{E}_{\sigma}, \text{dim}(T_x L \cap T_x \mathcal{E}_{\sigma})=n_1 \ .$$
As $T_x L$ has dimension $n_1+n_2$, it must contains $n_2$ vectors forming a basis of a supplement space for $T_x L \cap T_x \mathcal{E}_{\sigma}$ in $T_x L$. It means that $T_x L $ contains a collection of $n_2$ vectors, $$X_i=(X_{u,i},X_{q_1,i},X_{q_2,i},X_{p_1,i},X_{p_2,i})\in \R \times \R^{n_1} \times \R^{n_2} \times \R^{n_1} \times \R^{n_2}  \   ,$$  $\ i=1, \dots , n_2$, such that $\lbrace X_{q_2,i} \ , \ i=1, \dots , n_2 \rbrace $ is a basis of $\R^{n_2}$ . \\
The projection $\mathbf{p}_{\sigma}$ restricted to $L \cap \mathcal{E}_{\sigma}$ locally defines an immersion in $J^1(\R^{n_1},\R)$ if, for all $x=(u,q_1,q_2,p_1,p_2) \in L \cap \mathcal{E}_{\sigma}$, the following linear map is injective
\begin{align*}
 D_x(\mathbf{p}_{\sigma})_{\vert T_x(L \cap \mathcal{E}_{\sigma})} \ : & \ \qquad T_x(L \cap \mathcal{E}_{\sigma}) \qquad \longrightarrow \ T_{(u,q_1,p_1)}J^1(\R^{n_1},\R)\\
& (X_u,X_{q_1},0,X_{p_1},X_{p_2}) \longmapsto \ \ (X_u,X_{q_1},X_{p_1})
\end{align*} 

Its kernel consists in vectors $X=(0,0,0,0,X_{p_2})$ in  $T_x(L \cap \mathcal{E}_{\sigma})$. But, such a vector must satisfies the equation $(\star)$ together with any other vector of $T_x L$. In particular, considering the $n_2$ vectors $X_i$, $\ i=1, \dots , n_2$, it follows that the scalar product $X_{p_2}X_{q_2,i}$ must be $0$ for each $X_{q_2,i}$. Thus $X_{p_2}=0$, and we conclude that the space $\mathbf{p}_{\sigma} ( L \cap \mathcal{E}_{\sigma})$ in $J^1(\R^{n_1},\R)$ is locally the image of a submanifold by an immersion.\\

\noindent A similar proof holds for the contour case, exchanging the roles of $q_2$ and $p_2$.  \qed \\

\begin{rk}
A proof can also be made using generating function theory (section $3$). Moreover, we will use generating functions in section $3$ to prove the following.
\end{rk}

\begin{lem}\label{lem117}
For a generic Legendrian submanifold $L$, the slice $\sigma(L)$ of $L$ along $\R^{n_1}\times \lbrace 0 \rbrace$ (resp. contour $\kappa(L)$ in the direction of $\R ^{n_1}\times \lbrace 0 \rbrace$) is an immersed Legendrian submanifold. 
\end{lem}

\begin{prop}\label{prop11} \ Slice and contour are related by the following conjugation relations\footnote{Note that here again the same symbol $\T$ is used for the transformation $\T$ on different domains.}:
$$ \mathbf{T}\sigma(L)=\kappa(\mathbf{T}L) \qquad \text{ and } \qquad   \mathbf{T}\kappa(L)=\sigma(\mathbf{T}L) \ .$$
\end{prop}

\Pv \ It is a set theoretic check. \begin{align*}
\T \sigma (L) & = \T\lbrace (u,q_1,p_1) \ \vert \ \exists p_2 \in \R^{n_2}, (u,q_1,0,p_1,p_2) \in L \rbrace \\
& = \lbrace (q_1p_1-u,p_1,q_1) \ \vert \ \exists p_2 \in \R^{n_2}, (u,q_1,0,p_1,p_2) \in L \rbrace \\
& =  \lbrace (u,p_1,q_1) \ \vert \ \exists p_2 \in \R^{n_2}, (u,p_1,p_2,q_1,0) \in \T L \rbrace \\
& =  \kappa (\T L).
\end{align*}
The second equality follows from the first one and the fact that $\T$ is an involution. \qed \\

\vspace*{0.5cm}
\section{Sum and Convolution}\label{sec2}

In this section, we first define the operations sum and convolution of Legendrian submanifolds of $J^1(\R^n,\R)$ and give some of their properties. We then remind the operations on functions coming from convex analysis: the sum, the Legendre--Fenchel transform and the infimal convolution. We make precise the comparison between the operations on Legendrian submanifolds and those on functions lifted at the level of $1$-graphs. 

\subsection{Sum and Convolution of Legendrians submanifolds}\label{sec21}

\begin{df}\label{def21}
Let $L_1$ and $L_2$ be two Legendrian submanifolds of $J^1(\R^n,\R)$. The \textbf{sum} of $L_1$ and $L_2$ is the set $$L_1 \p L_2=\lbrace (u_1+u_2,q,p_1+p_2) \ \vert \ (u_i,q,p_i) \in L_i , i\in\lbrace 1,2 \rbrace\rbrace \ . $$
\end{df}

\begin{ex}\label{ex20}
Figure \ref{fig4} shows what may happen when making the sum of two opposite cusps. The middle case is the degenerate one. 
\begin{figure}[ht]
\begin{center}
 {\includegraphics[scale=0.55]{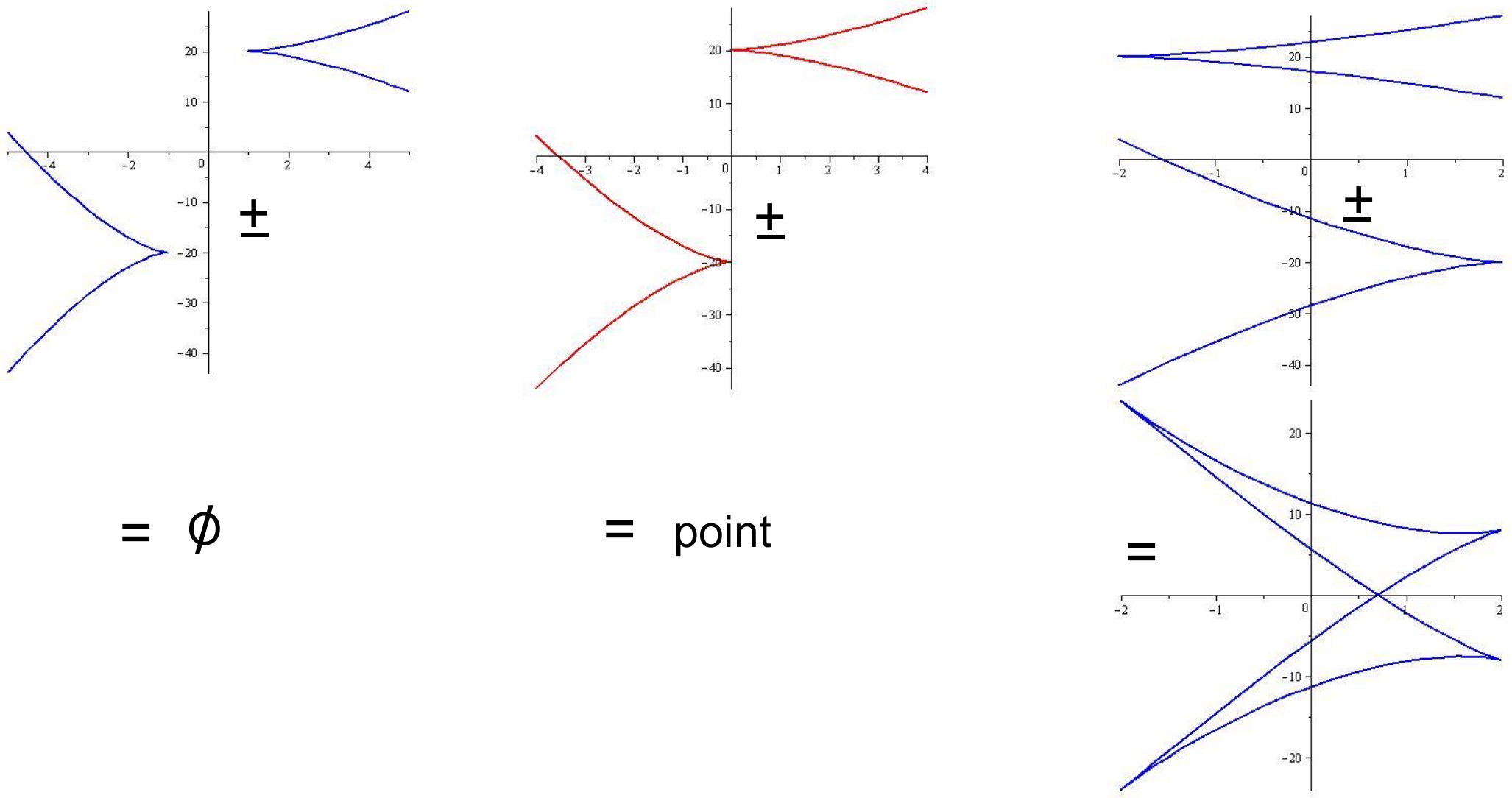}}
\end{center}
\caption{}
\label{fig4}
\end{figure}

\end{ex}

We will prove the following result using generating functions in section $3$.

\begin{lem}\label{lem24}
Let $L_1$ be a Legendrian submanifold of $J^1(\R^n,\R)$. For a generic Legendrian submanifold $L_2 \subset J^1(\R^n,\R)$, the sum $L_1 \p L_2$ is an immersed Legendrian submanifold.
\end{lem}

\begin{df}\label{def22}
Let $L_1$ and $L_2$ be two Legendrian submanifolds of $J^1(\R^n,\R)$. The \textbf{convolution} of $L_1$ and $L_2$ is the set $$L_1 \s L_2=\lbrace (u_1+u_2,q_1+q_2,p) \ \vert \ (u_i,q_i,p) \in L_i , i\in\lbrace 1,2 \rbrace\rbrace \ . $$
\end{df}

\begin{rk}\label{rk22}

Even when the sum produces embedded Legendrian submanifolds, non-generic singularities of wave fronts appear generically.

For example, in the case $n=2$, the sum of a line of cusp with another one such that their directions are transverse -- which is a generic case -- is a well embedded submanifold, but its wave front is systematically a \textsl{"handkerchief"} singularity, which is not generic, see Figure \ref{fig5}.

\begin{figure}[ht]
\begin{center}
 {\includegraphics[scale=0.55]{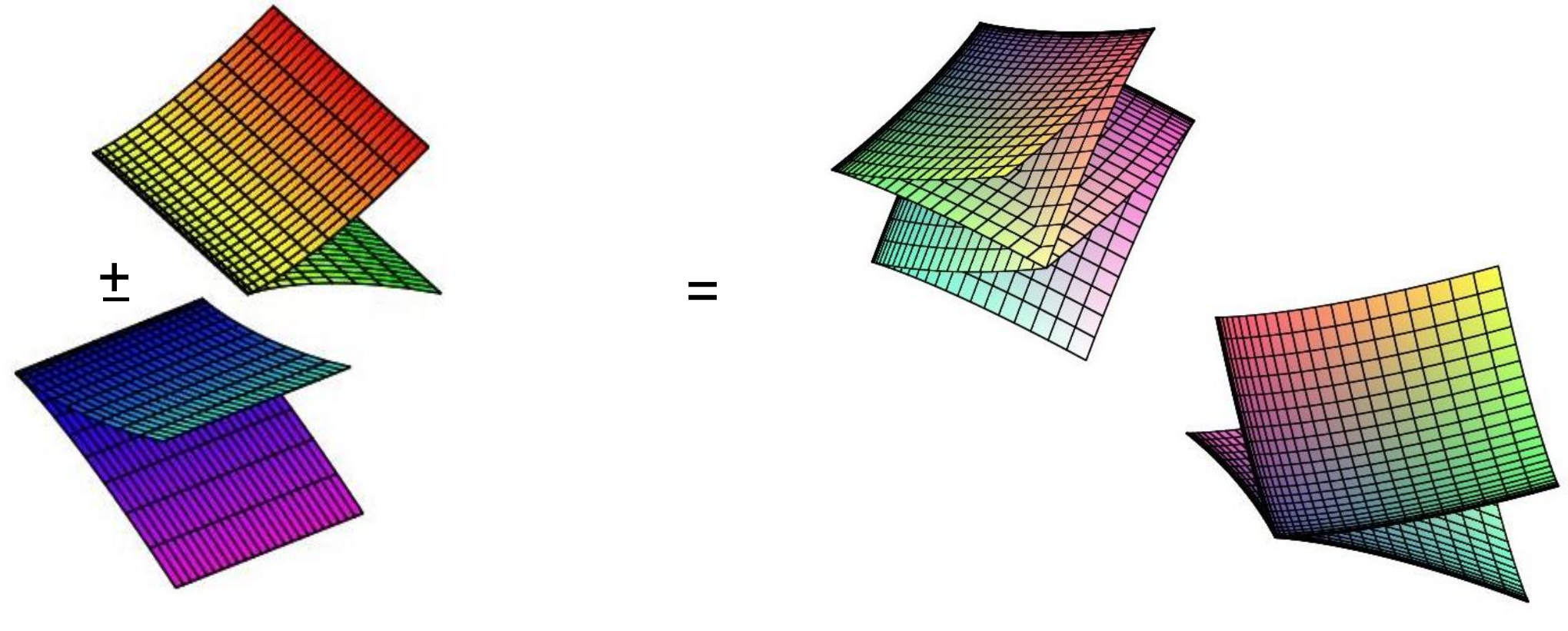}}
\end{center}
\caption{The sum of two lines of cusps produces a \textit{handkerchief}.} (the result is given with two different angles of view)
\label{fig5}
\end{figure}    
\end{rk}

\begin{rk}\label{rk21} The sum of $L_1$ and $L_2$ can be obtained by combining the operations product and slice as follow. Consider the image of the product $L_1 \ti L_2$ by the contactomorphism: \begin{align*}
\Phi \quad :  \quad J^1(\R^n \times \R^n , \R) \ & \longrightarrow \quad \quad J^1(\R^n \times \R^n , \R) \\
 (u,q_1,q_2,p_1,p_2) \ & \longmapsto  \ (u,\frac{q_1+q_2}{2},\frac{q_1-q_2}{2},p_1+p_2,p_1-p_2) \ ,
\end{align*}
$$\Phi(L_1\ti L_2)=\lbrace(u_1+u_2,\frac{q_1+q_2}{2},\frac{q_1-q_2}{2},p_1+p_2,p_1-p_2) \ \vert \ (u_i,q,p_i) \in L_i \rbrace  \ .$$ 

Then we compute the slice along $\R^n \times \lbrace 0 \rbrace$ of $\Phi(L_1\ti L_2)$: \begin{align*}
\mathbf{p}_{\sigma}(\Phi(L_1 \ti L_2) \cap \mathcal{E}_{\sigma}) & = \mathbf{p}_{\sigma}(\lbrace (u_1+u_2,q,0,p_1+p_2,p_1-p_2) \ \vert \ (u_i,q,p_i) \in L_i \rbrace \\
 & = \lbrace (u_1+u_2,q,p_1+p_2) \ \vert \ (u_i,q,p_i) \in L_i \rbrace \end{align*}
and recover the sum $L_1 \p L_2$.\\

The case of the convolution is similar: to obtain $ L_1 \s L_2 $ from $L_1 \ti L_2$ and the contour operation, we replace the contactomorphism $\Phi$ by: \begin{align*}
\Psi \quad :  \quad J^1(\R^n \times \R^n , \R) \ & \longrightarrow \quad \quad J^1(\R^n \times \R^n , \R) \\
 (u,q_1,q_2,p_1,p_2) \ & \longmapsto  \ (u,q_1+q_2,q_1-q_2,\frac{p_1+p_2}{2},\frac{p_1-p_2}{2}) \ ,
\end{align*} 
and take the contour of $\Psi(L_1 \ti L_2)$ in the direction of $\R^n \times \lbrace 0 \rbrace$ to get $L_1 \s L_2$.      
\end{rk}

\begin{thm}\label{thm21} Let $L_1$ and $L_2$ be two Legendrian submanifolds. Then $$\underline{i}) \ \mathbf{T}(L_1 \p L_2)=(\mathbf{T}L_1) \s (\mathbf{T}L_2) \quad \text{and} \quad \underline{ii}) \  \mathbf{T}(L_1 \s L_2)=(\mathbf{T}L_1) \p (\mathbf{T}L_2) .$$ 
\end{thm}

\Pv \ Foremost, note that $\Phi$ and $\Psi$ (Remark \ref{rk21}) are conjugated by $\T$:  \ $\Psi=\T \circ \Phi \circ \T \ .$ Then, using together Proposition \ref{prop11}, and Remarks \ref{rk12} and \ref{rk21}, we compute \begin{align*}
\underline{i)} \ \mathbf{T}(L_1 \p L_2) & = \T(\sigma(\Phi(L_1 \ti L_2))   \qquad \underline{ii)} \ \mathbf{T}(L_1 \s L_2) & = & \, \T(\kappa(\Psi(L_1 \ti L_2)) \\
& = \kappa(\T \circ \Phi(L_1 \ti L_2)) & = & \, \sigma(\T \circ \Psi(L_1 \ti L_2)) \\
& = \kappa(\Psi \circ \T (L_1 \ti L_2)) & = & \,  \sigma(\Phi \circ \T (L_1 \ti L_2)) \\
& = \kappa(\Psi ((\T L_1) \ti (\T L_2))   & = & \, \sigma(\Phi ((\T L_1) \ti (\T L_2))\\ 
& =(\mathbf{T}L_1) \s (\mathbf{T}L_2) & = & \, (\mathbf{T}L_1) \p (\mathbf{T}L_2)
\end{align*}
\qed \\

\begin{rk}
As a consequence, an analogous of Lemma \ref{lem24} holds for the convolution operation.
\end{rk}

\subsection{Legendre--Fenchel transform and infimal-convolution of functions}\label{sec22}
As $\nabla(f_1+f_2)=\nabla f_1 + \nabla f_2$, it is clear that the sum operation lifts the sum of functions at the level of the spaces of $1$-jets: $$(j^1f_1) \p (j^1f_2 )=j^1(f_1 + f_2).$$ 

\begin{df}\label{def23}
In that note, a function $f$ defined on $\R^n$ is said to be \textbf{simple} if:
$$\qquad \qquad \nabla f \text{ is a diffeomorphism of } \R^n .$$ 
\end{df}  

\begin{rk}
Note that a simple function admits a single critical point, which is non-degenerate in the sense of Morse theory.
\end{rk}

We remind then the essential facts about Legendre--Fenchel transform and infimal-convolution. The reader is referred to \cite{1} (chapters 12 and 13) and \cite{13} (chapter 2) for further details.

\begin{df} The \textbf{Legendre transform $f^{\mathbf{t}}$} of a simple function $f$ is defined by $$(1) \qquad \qquad  f^{\mathbf{t}}(q)=q(\nabla f)^{-1}(q) - f((\nabla f)^{-1}(q)) . \qquad \qquad$$
\end{df}

Formula (1) creates a function $f^{\mathbf{t}}$ whose gradient coincide with the inverse of the gradient of $f$: $$\nabla f^{\mathbf{t}}=(\nabla f)^{-1} .$$ 
In particular, $f^{\mathbf{t}}$ is also simple. Apply twice Legendre transform gives back the original function: $f^{\mathbf{t}\mathbf{t}}=f$.
Note that the widest parabolas are sent to the sharpest ones and reciprocally: $$q \mapsto \frac{a}{2}q^2
\text{ is sent onto } q \mapsto \frac{1}{2a}q^2 \ .$$

The following Lemma is easy to get.
\begin{lem}
Let $f$ be a simple function. Then $$\T(j^1 f)=j^1 (f^{\mathbf{t}}).$$
\end{lem}

Legendre--Fenchel transform extends the Legendre transform to functions which are \textit{almost} simple and convex.

\begin{df} 
Let $f$ be a function on $\R^n$. Its Legendre--Fenchel transform $f^*$ is defined by $$(2) \qquad \qquad f^*: \R^n \rightarrow \R \cup \lbrace + \infty \rbrace \ , \  q \in \R^n                                                                                                                                                                                                                                                                                                                                                                                                                        \mapsto \sup_{v \in \R^n}\lbrace qv-f(v)\rbrace . \qquad \qquad$$ 
\end{df}
 
Consider the $q$-parametrised family of functions $(F_q : v \mapsto qv-f(v))_{q \in \R^n}$. If $f$ is simple and convex, then each $F_q$ is simple and concave, and so admits a unique critical point $v_q$, whose value is a global maximum. Thus, the definitions by formulas (1) and (2) match perfectly when $f$ is simple and convex. 

\begin{figure}[ht]
\begin{center}
\includegraphics[scale=0.8]{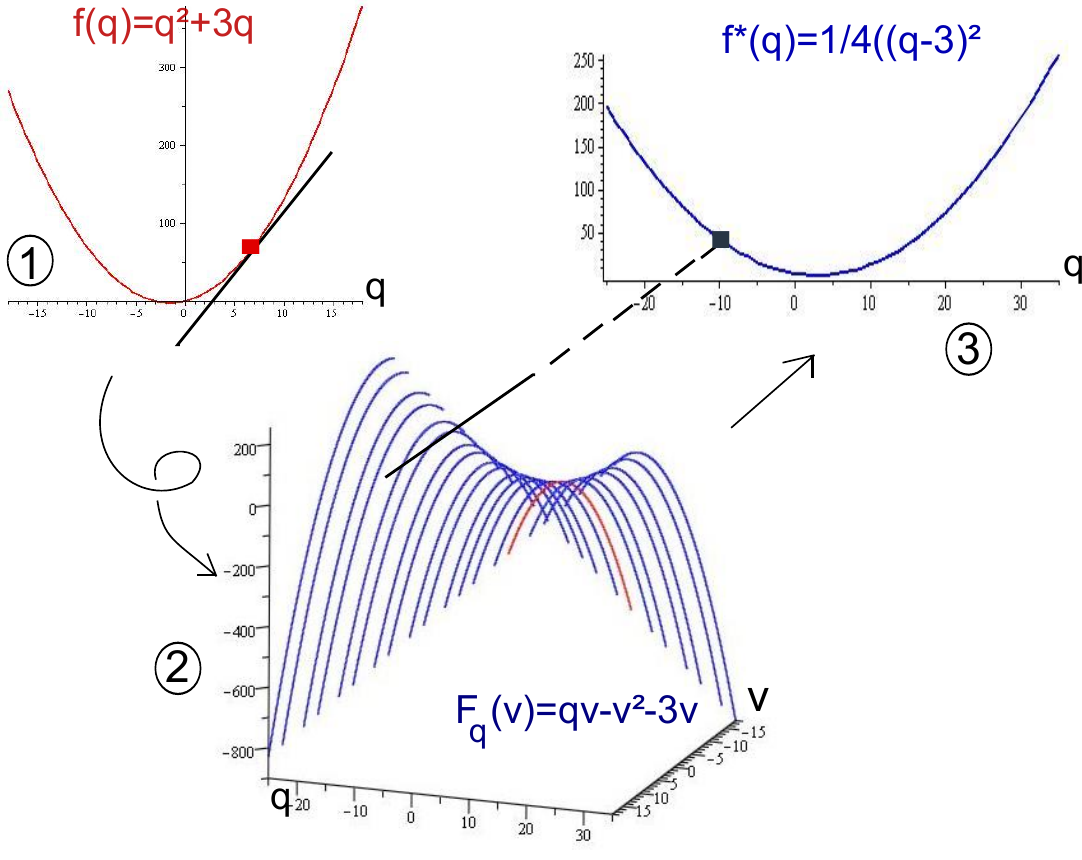}
\end{center}
\caption{Construction of the Legendre-Fenchel transform of $f : q \mapsto q^2+3q$.}
\label{fig6}
\end{figure}

The Legendre--Fenchel construction is illustrated by figure \ref{fig6} with a convex parabolic function defined on $\R$. First is represented the graph of $f: q \mapsto q^2 +3q$, then the $q$-parametrised family $(F_q)_{q \in \R^n}$, and finally, the graph of $f^*$ is obtained from the previous one by projecting the top of the concave parabola $v \mapsto qv-v^2-3v$ for each $q$.

From a geometrical point of view, observe that when $f$ is simple, the Legendrian submanifold $j^1f$ can also be viewed as a graph with respect to the $p$-variables. It is equivalent to observe that $ \T (j^1 f)$ is also the $1$-graph of a function. We may summerize this by the following

\begin{lem}
Let $f$ be a simple and convex function. Then $\T(j^1 f)$ is the $1$-graph of the Legendre--Fenchel transform of $f$ $$\T(j^1 f)=j^1 f^*.$$
Thus $f^*$ is also simple and convex.
\end{lem}

\begin{df} Let $f_1$ and $f_2$ be two functions. Their infimal-convolution $f_1 \square f_2$ is defined by $$f_1 \square f_2: \R^n \rightarrow \R \cup \lbrace - \infty \rbrace \ , \  q \mapsto \inf_{v \in \R^k}\lbrace f_1(v)+f_2(q-v) \rbrace. $$
\end{df}

If $f_1$ and $f_2$ are simple and convex, then for all $q \in \R^n$ the function $F_q : v \mapsto f_1(v)+f_2(q-v)$ is also a simple and convex function. So, for each $q$, $F_q$ admits a unique critical point $v_q$, whose value must be the absolute minimum of $F_q$: $$f_1 \square f_2 (q) =f_1(v_q)+f_2(q-v_q) \  \text{ for } v_q \text{ such that } \nabla f_1(v_q)=\nabla f_2(q-v_q) \ .$$ One obtains the definition of the convolution of the $1$-graphs $j^1f_1$ and $j^1f_2$ by rewriting $$f_1\square f_2(q)=f_1(q_1)+f_2(q_2) \ , \text{ for } q=q_1+q_2 \text{ and } \nabla f_1(q_1)=\nabla f_2(q_2) \ ,$$
and conclude the following

\begin{lem}\label{lem22}
Let $f_1$ and $f_2$ be two functions which are simple and convex. Their infimal-convolution $f_1 \square f_2$ is also a simple and convex function, and $$(j^1f_1)  \s (j^1f_2)=j^1(f_1 \square f_2).$$
\end{lem}

Thus, the geometric definitions of sum, convolution and transformation $\T$ that we introduced match the lifts of the convex analysis definition -- sum, infimal-convolution and Legendre--Fenchel transform -- when this \textit{strong convexity property} holds: $f$ is a convex function such that $\nabla f$ realises a diffeomorphism between the base and the space of slopes.

On the other hand, one can consider the operations sum, convolution and transformation $\T$ of $1$-graphs of functions and project on $J^0(\R^n,\R)$. Let us focus on the projection of the transformation $\T$ for the rest of this subsection. If $f$ is simple, then $\T (j^1f)$ remains a $1$-graph. So $\textbf{pr}\circ\T$ can be seen as an operation on functions if we restrict to simple ones.

If $f$ is not simple, the wave front $\textbf{pr}\circ \T(j^1 f)$ is not the graph of a function anymore, but a multi-valued graph (see Example \ref{ex24}, Figure \ref{fig8}). In section \ref{sec3} we will see how to recover single-valued functions with the notion of \textit{selector}. 

For now, we give three examples of Legendre--Fenchel transform for not simple and convex functions, and let us compare with $\textbf{pr}\circ \T$.

\begin{ex}\label{ex22}
Let $f$ be the exponential function: $f(q)=e^q \ , \ q\in \R$. It is smooth, and convex, but the space of slopes is reduced to $]0,+\infty[$.

On the one hand the computation of the classical Legendre--Fenchel transform gives
$$f^*(q)=\begin{cases}
\quad +\infty \hspace{1.15cm}\text{ , \ if } \ q<0\\
\qquad 0 \hspace{1.25cm} \text{ , \ if } \ q=0\\
q(\ln(q) -1) \ \ \text{ , \ if } \ q>0 
\end{cases}$$

On the other hand, the transformation $\mathbf{T}$ applied to the $1$-graph of $f$ gives as a result: $$ \mathbf{T}(j^1f)= \lbrace u=e^x(x-1),q=e^x,p=x) \ \vert \ x \in \R \rbrace \ , $$
which is also the $1$-graph of the function $q \mapsto q(\ln(q)-1)$, but only defined on $\R_{>0}$. 

So here, the difference between the graph of $f^*$ and the wave front of $\mathbf{T}(j^1f)$ resides in the way to deal with the unreached slopes. If one wants to rectify, it would be sufficient to consider $f^*$ restricted to its \textsl{domain} (i.e. the subset of $\R$ where it takes finite values).  
\end{ex}  

\begin{ex}\label{ex23} 
Let $f$ be the function defines on $\R$ as follow
$$f(q)= \begin{cases} 
(q+1)^2 \ , \ q \in \ ]-\infty ; -1]\\
\quad \ 0 \hspace{0.8cm} , \ q \in \ ]-1;1[\\
(q-1)^2 \ , \ q \in \ [1;+\infty[
\end{cases} $$  

It is a strongly degenerate case, with a whole interval where $f$ is constant, however it is of class $C^1$ on $\R$. In this example, the graph of $f^*$ and the wave front of $\mathbf{T}(j^1f)$ are the same. Notice that a singularity emerges when applying the Legendre--Fenchel transform or transformation $\mathbf{T}$ (see figure \ref{fig7}). Indeed, one can compute 
$$f^*(q)= \begin{cases} 
\frac{1}{4}q^2-q \ , \ q\leqslant 0\\
{}\\
\frac{1}{4}q^2+q \ , \ q\geqslant 0 
\end{cases} $$ 
 Moreover, if one applies one more time the Legendre-Fenchel transform, it leads to the original function $f$.
 
\begin{figure}[h]
\begin{center}
  \subfigure{\includegraphics[scale=0.18]{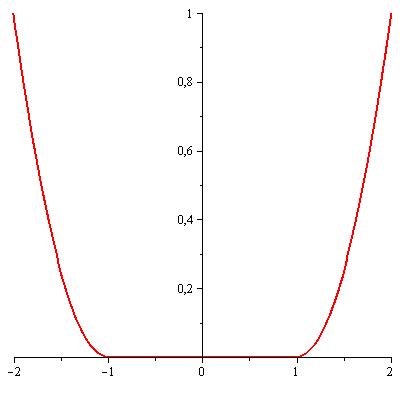}} \quad  $\overset{*}{\longrightarrow}$  \quad
  \subfigure{\includegraphics[scale=0.18]{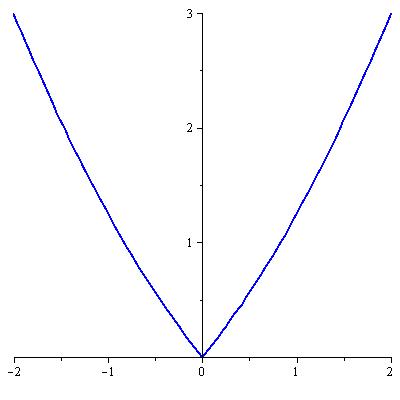}}\quad  $\overset{*}{\longrightarrow}$  \quad 
  \subfigure{\includegraphics[scale=0.18]{LF1}} 
\end{center}
\caption{successive transformations for a convex function with a constant part.}
\label{fig7}
\end{figure}
\end{ex}

\begin{ex}\label{ex24}
The fundamental difference is when we deal with \textsl{compact deformations} of simple and convex functions. Consider the function $f$ defined on $\R$ by  $$f(x)=q^4-3q^2 \ , $$ (see Figure \ref{fig1}), which coincides with a simple and convex function, except on a compact subset of $\R$ where it realises a concave bump. 

The definition by supremum of $f^*$ \textsl{selects} only the upper part of $\mathbf{pr}\circ \mathbf{T}(j^1f)$ to create the graph of $f^*$. Observe that a singularity appears, but this selection is continuous. Note also that the Legendre-Fenchel transform does not distinguish this case from the previous one.

If one applies one more time these operations -- see Figure \ref{fig8} -- on the one hand, using the transformation $\mathbf{T}$ permits to recover the full graph of $f$, while on the other hand, the Legendre-Fenchel transform gives as a result the convex hull of $f$.
\begin{figure}[h]
\begin{center}
\includegraphics[scale=0.6]{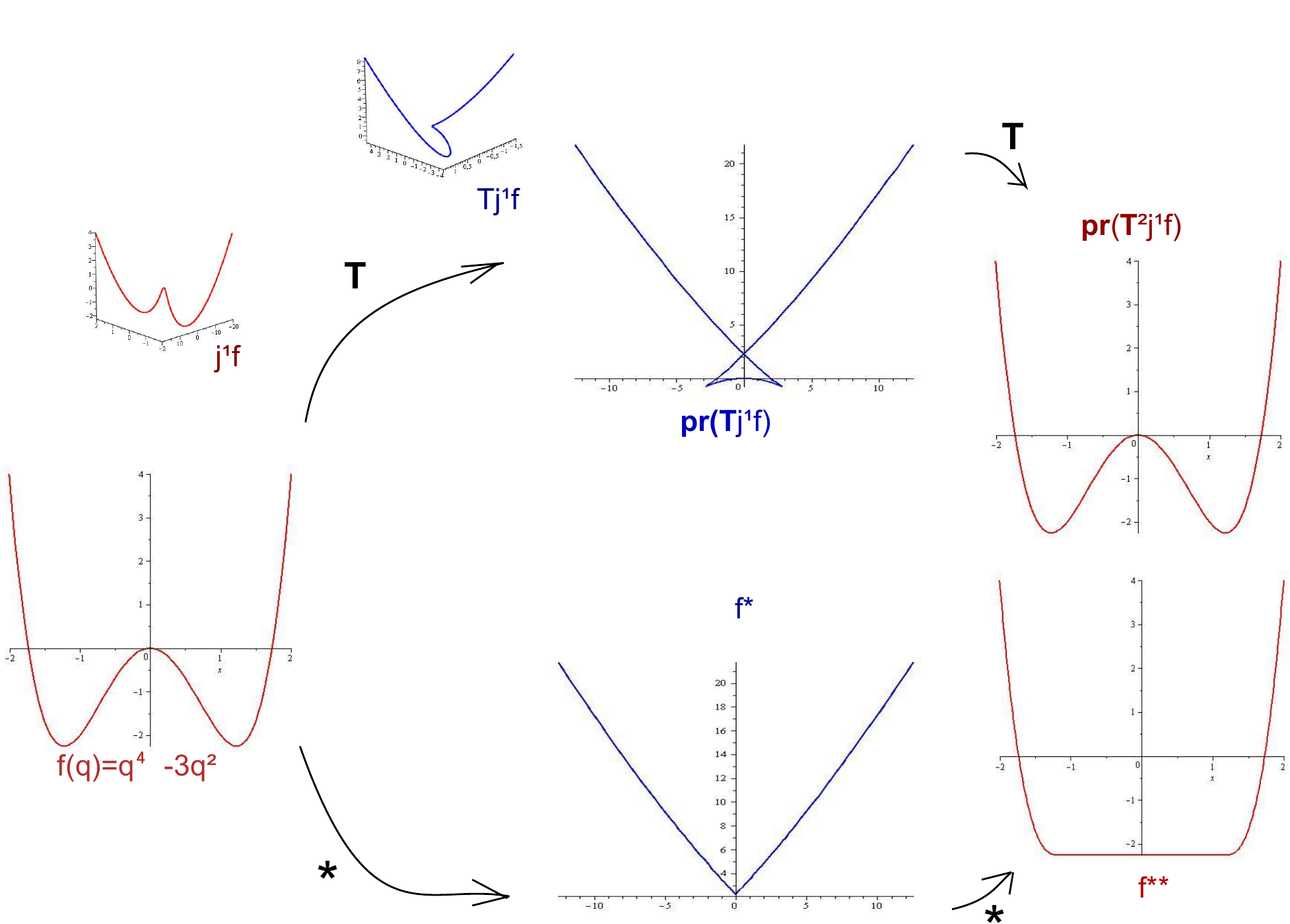}
\end{center}
\caption{}
\label{fig8}
\end{figure}
\end{ex}

\section{Generating Functions}\label{sec3}

\vspace{0.2cm}
\subsection{Generating Functions and Operations}

The $1$-graphs of functions defined on $\R^n$ are elementary examples of Legendrian submanifolds of $J^1(\R^n,\R)$ (Example \ref{ex11}). 
Generating functions define a larger class of Legendrian submanifolds, produced from $1$-graphs of functions thanks to a particular contour procedure. After reminding the classical notion of generating function, we see in this subsection that the set of Legendrian submanifolds that admit a generating function is stable by the three operations: transformation $\mathbf{T}$, sum and convolution (Lemma \ref{lem31}). \\

\noindent \textbf{Notations} \ Let $F$ be a function defined on a product $\R^n \times \R^k$ endowed with coordinates $(q,w)$. We denote by $\nabla_q F$ and $\nabla_w F$ the following maps : 
\begin{align*}
\nabla_q F : \ &\R^n \times \R^k  \longrightarrow \quad \R^n \quad & \nabla_w F : \ &\R^n \times \R^k \longrightarrow \quad \R^k\\
& (q,w)  \longmapsto \left( \frac{\partial F}{\partial q_j}(q,w)\right)_{j=1\ldots n} & & (q,w)  \longmapsto \left( \frac{\partial F}{\partial w_j}(q,w)\right)_{j=1\ldots k}
\end{align*}

\begin{df}\label{def31}
A \textbf{generating function} (\textbf{gf}) for a subset $L \subset J^1(\R^n,\R)$ is a function $F$ defined on $\R^n\times \R^k$,  $$q=(q_1,\ldots,q_n) \in \R^n  \text{ and } w=(w_1,\ldots,w_k) \in \R^k \ ,$$
such that $L$ is the contour (Definition \ref{def17}) along $\R^n \times \lbrace 0 \rbrace$ of $j^1F$. In that case, $L$ may be denoted by $L_F$, and it can be described as follow
$$L_F=\lbrace \left(u=F(q,w),q,p=\nabla_q F (q,w) \right) \ \vert \ \exists w \in \R^k \ , \ \nabla_w F (q,w)=0 \rbrace .$$
${}$\\
\noindent We say that $L=L_F$ is \textbf{the contour of F}. 
\end{df}

\noindent \textbf{Notation} \ For such an $F$, we bracket $\R^k$ so there will be no ambiguity between the base space and the additional parameter space, which disappears after the contour operation. Thus, a function $F$ defined on a product $\R^n \times \R^k$ may be understood as four different generating functions, defined respectively on $\R^n \times \R^k$, $\R^n \times (\R^k)$, $(\R^n) \times \R^k$ or $(\R^n \times \R^k)$  \\

\begin{rk}\label{rk32}${}$\\
$\bullet$ The $1$-graphs of functions are elementary examples, with $k=0$.\\
$\bullet$ Conversely, a Legendrian submanifold can always be obtained locally by such a contour procedure \cite{2}. However the existence of a generating function for a whole Legendrian submanifold is a strong global constraint. 
\end{rk}

\begin{df}\label{def32} 
For $F$ a generating function, with additional variable $w \in \R^k$, we will say that $F$ \textbf{satisfies the condition  $(\ast)$} if : $$(\ast) \quad \quad \nabla_w F \text{ is transverse to } \lbrace 0 \rbrace \ . \quad \quad$$ 
\end{df}

\begin{rk}\label{rk31}
Note that the condition $(\ast)$ corresponds to the transversality condition in Lemma \ref{lem11}. Thus, as soon as the condition $(\ast)$ is satisfied, $L_F$ is an immersed Legendrian submanifold of $J^1(\R^n,\R)$.\\  
\end{rk}

The following fact is folklore: 

\begin{prop}
The $(\ast)$ condition is generic for generating functions. 
\end{prop}

\noindent \textit{Idea of the proof.} This is a consequence of Thom's lemma for transversality in spaces of jets, \cite{19}. The condition $(\ast)$ can be reformulated asking that the $2$-jet of $F$ does not meet some stratified submanifold $\mathcal{I}$ of the space of $2$-jets $J^2(\R^n \times \R^k,\R)$. Thom's lemma states that the $2$-jet of a generic function $F$ defined on $\R^n \times \R^k$ meets any submanifold of the space of $2$-jets $J^2(\R^n \times \R^k,\R)$ transversally. Counting that each stratum of $\mathcal{I}$ has codimension greater than the dimension of $j^2 F$ leads to conclude that the transversality means the absence of intersection. The number of strata of $\mathcal{I}$ is finite, which implies that, for a generic $F$, $j^2 F$ meets $\mathcal{I}$ transversally. \qed  \footnote{The same kind of proof -- combined with Lemma \ref{lem11} and the fact that every Legendrian submanifold can be locally realized by a generating function -- permits to show that the two reductions operations, slice and contour (section 1), produce generically immersed Legendrian submanifolds when transversality holds at infinity.}\\

The operations considered so far behave well with respect to the notion of generating function.\\

\begin{lem}\label{lem31} $(1)$ \ If a Legendrian submanifold $L$ of $J^1(\R^n,\R)$ admits $F$ defined on $\R^n \times (\R^k)$ as a gf, then $\mathbf{T}L$ is the contour of \begin{align*}
F_{\mathbf{T}} : \ & \R^n \times (\R^n \times \R^k)  \longrightarrow  \ \R \\ 
 & \ \ (q,v,w) \ \longmapsto  qv-F(v,w)
\end{align*}
defined on $\R^n \times (\R^n \times \R^k)$. \\

\noindent $(2)$ \ If $L$ is a Legendrian submanifold of $J^1(\R^{n_1}\times\R^{n_2},\R)$ is the contour of $F$ defined on $\R^{n_1} \times \R^{n_2} \times (\R^k)$, then\\  
- \ the slice $\sigma(L)$ along $\R^{n_1} \times \lbrace 0\rbrace$ is the contour $F_{\sigma}$ defined on $\R^{n_1} \times (\R^k)$ by: $$F_{\sigma}(q_1,w)=F(q_1,0,w) ,$$ \\
- \ the contour $\kappa(L)$ in the direction of $\R^{n_1} \times \lbrace 0\rbrace$ is the contour of $F_{\kappa}$ defined on $\R^{n_1} \times (\R^{n_2} \times \R^k)$ by $$F_{\kappa}(q_1,v,w)=F(q_1,v,w).$$

\noindent $(3)$ \ If $L_1$ and $L_2$ are Legendrian submanifolds of respectively $J^1(\R^{n_1},\R)$ and $J^1(\R^{n_2},\R)$, such that $L_1$ is the contour of $F_1$ defined on $\R^{n_1} \times (\R^{k_1})$, and $L_2$ is the contour of $F_2$ defined on $\R^{n_2} \times (\R^{k_2})$, then the product $L_1 \ti L_2$ is the contour of $F_{1 \ti 2}$ defined on $\R^{n_1} \times \R^{n_2} \times ( \R^{k_1} \times \R^{k_2} )$ by 
$$F_{1 \ti 2}(q_1,q_2,w_1,w_2)=F_1(q_1,w_1)+F_2(q_2,w_2). $$

\noindent $(4)$ \ If $L_1$ and $L_2$ are Legendrian submanifolds of $J^1(\R^n,\R)$ which are the contours of respectively $F_1$ defined on $\R^{n} \times(\R^{k_1})$, and $F_2$ defined on $\R^n \times (\R^{k_2})$, then \\
- \ the sum $L_1 \p L_2$ is the contour of $F_{1 \p 2}$ defined on $\R^n \times (\R^{k_1} \times \R^{k_2})$ by $$F_{1 \p 2}(q,w_1,w_2)=F_1(q,w_1)+F_2(q,w_2) ,$$
- \ the convolution $L_1 \s L_2$ is the contour of $F_{1 \s 2}$ defined on $\R^n \times (\R^n \times \R^{k_1} \times \R^{k_2})$ by $$F_{1 \s 2}(q,v,w_1,w_2)=F_1(v,w_1)+F_2(q-v,w_2).$$
\end{lem}

\Pv \ It is an elementary check. \qed

\begin{rk}\label{rk33}
If $L$, $L_1$ and $L_2$ are respectively the $1$-graphs of functions $f$, $f_1$ and $f_2$ defined on $\R^n$, then $F=f$, $F_1=f_1$ and $F_2=f_2$ are generating functions for respectively $L$, $L_1$ and $L_2$. Note that the above generating functions $F_{\mathbf{T}}$ and $F_{1\s 2}$ are the expressions which appear in the definitions of respectively the classical Legendre-Fenchel transform and the infimal-convolution.\\ 
\end{rk}

Let us remind Lemma \ref{lem117} and Lemma \ref{lem24}. We have at our disposal the local tool of generating functions to prove them.\\

\noindent \textit{Proof of Lemma \ref{lem117}.} \ $L$ can be locally expressed as the contour of a generating function $F$ defined on a product $\R^{n_1 +n_2}\times (\R^k)$. The (local) transversality condition may be translated as the transversality of $j^1 F$ to a sub-space of $J^1(\R^{n_1 +n_2}\times \R^k, \R)$. Thom's lemma ensures that this transversality condition is generically true. Then, reasoning with a countable covering of $J^1(\R^{n_1+n_2} \times \R^k,\R)$ by compact sets permits to conclude that the transversality is generic for the whole Legendrian submanifold. \qed \\

\noindent \textit{Idea of proof for Lemma \ref{lem24}.} \ As in the proof of Lemma \ref{lem117}, we express locally the two Legendrian submanifolds $L_1$ and $L_2$ as the contour of generating functions $F_1$ and $F_2$. So the genericity issue is again removed to a space of smooth applications. Fixing $F_1$, and using the Thom's method of \textit{"transversality under constraint"} presented by Laudenbach in \cite{7.1} section 5, one can show that a generic function $F_2 \in C^{\infty}(\R^n \times \R^{k_1})$ which satisfies the $(*)$ condition in $\R^{k_2}$ is such that $F_{1 \p 2}$ satisfies the $(*)$ condition in $\R^{k_1} \times \R^{k_2}$. The result follows again by the existence of a countable cover of $J^1(\R^n,\R)$ by compact sets. \qed \\

\begin{rk}\label{rk34} $\bullet$ Using successively the formulas given in Lemma \ref{lem31}, we deduce on the one hand a generating function for the \textit{contour of the slice}, and on the other hand a generating function for the {slice of the contour}. Observe that they match perfectly.

\[
\xymatrixrowsep{0cm}
\xymatrixcolsep{0cm}
\xymatrix{
 & F_{\sigma}(q_1,q_2,w) {} \qquad \quad & &\\
 & =F(q_1,q_2,0,w) {} \qquad \quad \ar@/^2pc/[rd]^{\kappa} & &\\
   & & (F_{\sigma})_{\kappa}\\
   F(q_1,q_2,q_3,w) {} \qquad \quad \ar@/^2pc/[ruu]^{\sigma}\ar@/_2pc/[rdd]^{\kappa} &  & {} \qquad \ =F(q_1,v,0,w) & \textcolor{forestgreen}{\surd} \ar@/^2pc/@*{[forestgreen]}[ld] \ar@/_2pc/@*{[forestgreen]}[lu]\\
  & &  (F_{\kappa})_{\sigma} & \\
& F_{\kappa}(q_1,v,q_3,w) {} \qquad \quad  \ar@/_2pc/[ru]^{\sigma} & &\\
& =F(q_1,v,q_3,w) {} \qquad \quad  & &
}
\]

${}$

In other words, the slice and the contour commute at the level of generating functions.\\

$\bullet$ In contrast, let us compare the following generating functions:
\[
\xymatrixrowsep{0cm}
\xymatrixcolsep{0cm}
\xymatrix{
 1.) \qquad & F_{\mathbf{T}}(q_1,q_2,v_1,v_2,w) \quad \ar@/^2pc/[rdd]^{\sigma} & &\\
 & =q_1v_1+q_2v_2-F(v_1,v_2,w)  & & \textcolor{forestgreen}{\surd} \ar@/^2pc/@*{[forestgreen]}[ldddd] \ar@/_2pc/@*{[forestgreen]}[ld]\\
   & & (F_{\mathbf{T}})_{\sigma}(q_1,v_1,v_2,w) & \\
   F(q_1,q_2,w) {} \qquad \quad \ar@/^2pc/[ruuu]^{\mathbf{T}}\ar@/_2pc/[rddd]^{\kappa} & & {} \qquad \qquad =q_1v_1-F(v_1,v_2,w) & \\
  & &  (F_{\kappa})_{\mathbf{T}}(q_1,V,v,w) & \\
& F_{\kappa}(q_1,v,w) {} \qquad \quad   & \qquad \qquad =q_1V-F(V,v_1,w) & \\
& =F(q_1,v,w) \quad \ar@/_2pc/[ru]^{\mathbf{T}} &  &
}
\]\\

\[
\xymatrixrowsep{0cm}
\xymatrixcolsep{0cm}
\xymatrix{
 2.) \qquad & F_{\mathbf{T}}(q_1,q_2,v_1,v_2,w) \quad \ar@/^2pc/[rdd]^{\kappa} & &\\
 & =q_1v_1+q_2v_2-F(v_1,v_2,w)  & & \\
   & & (F_{\mathbf{T}})_{\kappa}(q_1,V,v_1,v_2,w) & \\
   F(q_1,q_2,w) {} \qquad \quad \ar@/^2pc/[ruuu]^{\mathbf{T}}\ar@/_2pc/[rddd]^{\sigma} & & {} \qquad \qquad =q_1v_1+Vv_2-F(v_1,v_2,w) & \\
  & & (F_{\sigma})_{\mathbf{T}}(q_1,v,w) \qquad & \\
& F_{\sigma}(q_1,w) {} \qquad \quad   &  \qquad =q_1v-F(v,0,w) & \textcolor{red}{\times}  \ar@/^2pc/@*{[red]}[l] \ar@/_2pc/@*{[red]}[luu]\\
& =F(q_1,0,w) \quad \ar@/_2pc/[ru]^{\mathbf{T}} & &
}
\]

In case 1.) the resulting generating functions are the same, whereas in case 2.) the two generating functions  significantly differ.\\

The same dissymmetry can be observe dealing with generating functions for the sum and the convolution. Starting with $F_1(q,w_1)$ and $F_2(q,w_2)$, compare \\

\noindent $\underline{i}) \ F_{(\mathbf{T}1) \p (\mathbf{T}2)} \text{ with } F_{\mathbf{T}(1 \s 2)}: $\\
\begin{align*}
& F_{(\mathbf{T}1) \p (\mathbf{T}2)}(q,v_1,v_2,w_1,w_2) & = \ & \ q(v_1+v_2)-(F_1(v_1,w_1)+F_2(v_2,w_2)) \qquad\\
& F_{\mathbf{T}(1 \s 2)}(q,V,v,w_1,w_2) & = \ & \ qV-(F_1(v,w_1)+F_2(V-v,w_2))  
\end{align*} \begin{flushright}
$\textcolor{forestgreen}{\surd}$\\
\end{flushright}

\vspace{0.2cm}
\noindent $\underline{ii}) \ F_{(\mathbf{T}1) \s (\mathbf{T}2)} \text{ with } F_{\mathbf{T}(1 \p 2)}: $ 

\begin{align*}
& F_{(\mathbf{T}1) \s (\mathbf{T}2)} (q,V,v_1,v_2,w_1,w_2) & = \ & \ V(v_1-v_2)+qv_2 -(F_1(v_1,w_1)+F_2(v_2,w_2)) \\
& F_{\mathbf{T}(1 \p 2)} (q,v,w_1,w_2) & = \ & \ qv-(F_1(v,w_1)+F_2(v,w_2)) 
\end{align*}
\begin{flushright}
$\textcolor{red}{\times}$\\
\end{flushright}

We come back to these comparisons at the end of the paper.\\

Note that even if the formulas differ, a direct computation in each case ensures that the two generating functions realise the same subsets. In particular, it leads to another proof of the Theorem \ref{thm21}.
\end{rk}

\vspace{0.2cm}
\subsection{Selector and operations}\label{sec32} We first recall the notion of min-max selector. 

\subsubsection{Construction}
The \textit{selector}, or \textit{min-max selector} admits various frameworks and uses \cite{5}, \cite{3} and \cite{12}. When working with functions defined on non-compact sets, some constraint has to be added to handle the Morse dynamic at infinity. In this note we have chosen a standard model at infinity based on \textit{simple functions}, rather than the classical quadratic one. Not only is this notion more general, but it also arises naturally working with the Legendre transform and the transformation $\T$. 

\subsection*{Min-Max for functions almost-simple} Let $f$ be a function defined on $\R^k$ such that its gradient map $\nabla f$ is a proper map. It satisfies a strong version of the Palais-Smale condition $$\text{if there exists a sequence } (w_n)_n \text{ in } \R^k \text{ such that } \nabla f (w_n) \underset{n \rightarrow + \infty}{\longrightarrow}0 \ ,$$ 
$$\text{then there exists a convergent subsequence } w_{\phi(n)}\underset{n \rightarrow + \infty}{\longrightarrow}\zeta \ .$$ Moreover, the limit $\zeta$ must be a critical point. As a consequence, the set of critical points is compact, and so is the set of critical values.\\

\noindent Such a function $f$ permits to work with Morse theory, in the sense that: \\
\noindent 1. \ the topology of the sublevel sets of $f$ changes exactly when passing a critical value (Palais-Smale condition),\\
\noindent 2. \ these changes are gathered in a compact set of $\R^k$.\footnote{Remind that excellent Morse functions defined on a compact set form a $C^0$-dense subset of the space of functions \cite{9}.}\\

\noindent \textbf{Notations} \ For $\lambda \in \R$, we denote by $f^{\lambda}$ the sublevel set of $f$,
$$f^{\lambda}=f^{-1}(]-\infty,\lambda]) \ .$$
For two real numbers $\lambda_1\leqslant \lambda_2$, we denote the relative sublevel set by $$f_{\lambda_1}^{\lambda_2}=f^{-1}([\lambda_1,\lambda_2]) \ .$$

\begin{df} Let $f$ be a function such that $\nabla f$ is a proper map. Let $C$ denote a positive constant such that the set of critical values of $f$ is contained in $]-C,C[$. Then, $\forall \lambda >C $, we denote the sublevels sets $f^{\lambda}$ and $f^{-\lambda}$ by $f^{+\infty}$
and $f^{-\infty}$ respectively. We call $f^{+\infty}$ the \textbf{high sub-level} of $f$, resp.
$f^{-\infty}$ the \textbf{low sub-level} of $f$.
\end{df}

\begin{rk}
These notations and definitions give us a convenient way to work with relative homologies $H_*(f^{+\infty},f^{-\infty};\Q)$ and $H_*(f^{\lambda},f^{-\infty};\Q)$ when $\lambda$ varies in $[-C,C]$.
\end{rk}

If the relative homology groups $H_*(f^{+\infty},f^{-\infty};\Q)$ are all zero except in one degree $\iota \in \llbracket 0,k \rrbracket$ where it is equal to $\Q$, then one can look for the smallest $\lambda$ where a generator for $H_{\iota}(f^{+\infty},f^{-\infty};\Q)$ appears in $f^{\lambda}$.

\begin{df}
Let $f$ be a function such that $\nabla f$ is a proper map. If there is $\iota \in \llbracket 1,k \rrbracket $ such that $$H_*(f^{+\infty},f^{-\infty};\mathbb{Q}) \simeq \begin{cases}
\ \mathbb{Q} \quad \text{ if } * = \iota \\
\ \, 0 \, \quad \text{ if } * \neq \iota 
\end{cases}$$ we will say that $f$ has \textbf{simple relative homology of index $\iota$}.
In that case, we define \textbf{the min-max of $f$} as the real $$s(f):=\inf\lbrace \lambda  \ \vert \  H_{\iota}(f^{\lambda},f^{-\infty}) \overset{\mathbf{i}_*^{\lambda}}{\rightarrow }H_{\iota}(f^{+\infty},f^{-\infty}) \text{ is surjective}\rbrace \ , $$
where $\mathbf{i}_*^{\lambda}$ is the map induced by the inclusion map $$(f^{\lambda},f^{-\infty}) \overset{\mathbf{i}^{\lambda}}{\hookrightarrow}(f^{+\infty},f^{-\infty}) \ .$$
\end{df}

\begin{rk}${}$\\
\noindent $\bullet$ \ If $\lambda$ is not a critical value of $f$, the sublevel set $f^{\lambda}$ can be deformed by retraction onto the sublevel set $f^{\lambda - \epsilon}$ until $\lambda - \epsilon$ takes the value of a critical value $c$ of $f$. So the min-max is necessarily a critical value. \\ 
\noindent $\bullet$ \ Moreover, a relative cycle in $f^{s(f)}$ which realises a generator for $H_{\iota}(f^{+\infty},f^{-\infty})$ must pass through a critical point $\zeta \in f^{-1}(\lbrace s(f) \rbrace)$ which has Morse index $\iota$.
\end{rk}

\begin{df}\label{def35} Let $g$ be a simple function (Definition \ref{def23}). The \textbf{Morse index of g} is the Morse index of its critical point.
\end{df}

\begin{df}\label{def36}
A function $f$ defined on $\R^k$ is \textbf{almost-simple} if it can be decomposed into a sum as follows $$f=g+h \text{ such that } \begin{cases}g \text{ simple function }\\
 \exists B>0 \ \vert \ \forall w \in \R^k \ , \ \parallel \nabla h(w) \parallel \leqslant B. \end{cases}$$ 
\end{df}

\begin{ex}\label{ex31}
In particular, non-degenerate quadratic forms are (fundamental) examples of simple functions, with critical value zero and Morse index equals to the index of the quadratic form. Functions which are quadratic at infinity -- i.e. equal to a non-degenerate quadratic form outside of a compact set (as consider in \cite{3}) -- are first examples of almost-simple functions.  
\end{ex}

\begin{rk}${}$\\
$\bullet$ A simple or almost-simple function has a proper gradient map.\\
$\bullet$ A simple function of Morse index $\iota$ has simple relative homology of index $\iota$.
\end{rk}

When two functions $f$ and $g$, with proper gradients, are such that $\nabla(f-g)$ is bounded, the Moser's path method permits to construct a diffeomorphism which maps the big (non-critical) sublevel sets of $f$ onto those of $g$ -- see for instance the proofs of Proposition 11 in \cite{14} and Lemma C.2 in \cite{12}. 

\begin{lem}
If $f$ and $g$ are two functions such that $\nabla f$ and $\nabla g$ are proper maps, and the difference $\nabla (f-g)$ is bounded, then their high (respectively low) sublevel sets are diffeomorphic.
\end{lem}

\begin{cor}
If $f$ is a almost-simple function which can be decomposed into $g+h$ as in Definition \ref{def36}, with $g$ a simple function of Morse index $\iota$, then $f$ has a simple relative homology of index $\iota$. Thus: \\

\noindent 1. \ $f$ admits a min-max, denoted $s(f)$\\

\noindent 2. \ if $f=\tilde{g}+\tilde{h}$ is another decomposition of $f$ as in Definition \ref{def36}, then the Morse index of the simple function $\tilde{g}$ must be $\iota$. We name $\iota$ \textbf{the Morse index of the almost-simple function $f$}.
\end{cor}

\begin{df}\label{def37}
A function $f$ defined on $\R^k$ is \textbf{almost-convex} (respectively \textbf{almost-concave}) if it is almost-simple of Morse index 0 (resp. Morse index $k$). 
\end{df}

The following is part of folklore.

\begin{lem}\label{lem33}
If $f$ is a almost-convex function (respectively almost-concave) then the min-max of $f$ is its absolute minimum (resp. absolute maximum): $$s(f)=\min f $$ $$(\text{resp. } \  s(f)=\max f ) \ . \qquad  $$
\end{lem}

\Pv \ Suppose $f$ is almost-convex. Note that $f^{-\infty}$ is empty, so the relative homology $H_*(f^{+\infty},f^{-\infty}; \Q)$ is the same as $H_*(f^{+\infty}; \Q)$, and that the high sublevel set is a disk. The minimum $m$ of $f$ is obviously the smallest candidate among the critical values of index $0$. Moreover, passing the level $m$ gives rise to the sublevel set $f^{m+\epsilon}$, which is a disk included in the disk $f^{+\infty}$ for $\epsilon$ small enough, so a deformation retract of the high sublevel set. Thus, $H_{0}(f^{m+\epsilon},f^{-\infty}) \overset{\mathbf{i}_*^{m+\epsilon}}{\rightarrow }H_{0}(f^{+\infty},f^{-\infty})$ is bijective.\\

\noindent If $f$ is almost-concave, $f_{-\infty}^{+\infty}$ is a disk attached to the low sublevel set $f^{-\infty}$. Denote by $M$ the maximum of $f$. Suppose there exists a critical value $c_1<M$ such that $f^{c_1}$ contains a relative cycle generating $H_k(f^{+\infty},f^{-\infty}; \Q)$. Then there exists $\epsilon>0$ such that $c_1+\epsilon<M$ and $f^{c_1+\epsilon}$ is non critical. The relative sublevel set $f^{c_1+\epsilon}_{-\infty}$ should have the same homology as the disk $f^{+\infty}_{-\infty}$. But, since $c_1+\epsilon$ is note the maximum, the level set $f^{-1}(c_1+\epsilon)$ is non-empty, and $f^{c_1 +\epsilon}_{-\infty}$ has more than one boundary component. It is absurd. So the smallest critical sublevel set where a generator for $H_k(f^{+\infty},f^{-\infty};\Q)$ appears is $f^M$. \qed \\

\subsection*{Selector for generating functions almost-simple} \ This is the $q$-parametrized version of the min-max.

\begin{df}
A generating function $G$ defined on $\R^n \times \R^k$ is \textbf{simple} if:
$$\forall q \in \R^n \ , \ G(q,.) \text{ is a simple function} \ .$$
\end{df}  

One can prove the following as a consequence of a \textit{parametrized version of the Morse Lemma}, \cite{27}.

\begin{lem}\label{lem34}
If $G$ is a simple generating function, the index of $G(q,.)$ must be the same for all $q \in \R^n$. 
\end{lem}

\begin{df}\label{def39}
A generating function $F$ defined on $\R^n \times \R^k$ is said \textbf{almost-simple} (gf \textbf{AS}) if it can be decomposed as follow: 
$$ F=G+H \text{ , with } \forall q \in \R^n \  \begin{cases}G(q,.) \text{ simple }\\
\parallel \nabla_w H(q,.) \parallel \leqslant B(q)
\end{cases} \ ,$$
where $B$ is a continuous bounded function defined on $\R^n$.
\end{df}

\begin{cor}\label{cor31} ${}$\\
1. If $F$ is a gf AS, then $F(q,.)$ is a almost-simple function for all $q\in \R^n$.\\
2. If $F$ is a gf AS, the Morse index of $F(q,.)$ is the same for all $q \in \R^n$. 
\end{cor}

\Pv \ The first point comes directly from Definition \ref{def39}. The second follows from the Lemma \ref{lem34}.\qed  \\

\begin{df}
Let $F$ be a gf AS defined on $\R^n \times \R^k$. For all $q \in \R^n$, denote by $s(F)(q)$ the min-max of the almost-simple function $F(q,.)$: $$s(F)(q):= s(F(q,.)) \ .$$ The \textbf{selector of $F$} is \begin{align*}
s(F) \ : \ & \R^n \longrightarrow \ \R \\
& \ q \ \longmapsto s(F)(q)
\end{align*}
\end{df}

The following Lemma is also part of the folklore. One must use the inequality in Definition \ref{def39} to reduce the issue in our framework to the compact case proved in \cite{7.1}, section 6.3. 

\begin{lem}
The selector of a gf AS is continuous.
\end{lem}

%\Pv \ If $f$ and $g$ are two smooth functions defined on a compact set such that the min-max's $s(f)$ and $s(g)$ exist, then one can show that\footnote{See \cite{7.1}, section 6.3.} $$\forall \epsilon > 0, \exists \eta
%>0 \text{ s.t. } \parallel f-g \parallel_{C^0}<\eta \Rightarrow \vert s(f)-s(g) \vert <\epsilon.$$
%
%\noindent Thus, it is sufficient to show that the continuity issue for the selector of a gf AS can be reduced in a compact set. Let $F$ be a gf AS defined on $\R^n \times (\R^k)$. Given $q_0\in \R^n$ and $r>0$, consider the closed ball $\bar{B}(q_0,r) \subset \R^n$. We show that the union over $q \in \bar{B}(q_0,r)$ of all critical points of the functions $F(q, \cdot )$ is stuck in a compact set $K$. Then, for $q \in \bar{B}(q_0,r)$, we restrict the functions $F(q,\cdot)$ to $K$ and apply the above argument using the $C^0$-topology based on $K$. 
%
%If $(q,w) \in \bar{B}(q_0,r)\times \R^k$ satisfies $\nabla_w F(q,w)=0$, then $$\parallel \nabla_w G(q,w) \parallel \leqslant B(q) ,$$ and thus $w \in \left( \nabla_w G(q,\cdot )\right)^{-1}\left(\bar{B}(0,B(q))  \right) \subset \left( \nabla_w G(q,\cdot )\right)^{-1}\left(\bar{B}(0,M)  \right)$ where $M$ denotes the maximum of the continuous bounded function $B$ on $\bar{B}(q_0,r).$ 
%The map \begin{align*}
%\bar{B}(q_0,r) \times \R^k & \longrightarrow \ \R^k \\
%(q,W) \ & \longmapsto \left( \nabla_w G (q,\cdot )\right)^{-1}(W)
%\end{align*} 
%is continuous, thus the image of $\bar{B}(q_0,r) \times \bar{B}(0,M))$ is a compact subset $K \subset \R^k$. \qed

\vspace{0.5cm}
\subsubsection{Selector and sum}

\begin{prop}\label{prop31}
Let $F_1$ and $F_2$ be two gf AS. Then $F_{1 \p 2}$ is a gf AS and: $$ s(F_{1 \p 2})=s(F_1)+s(F_2)  \ .$$  
\end{prop}

\Pv \ It can be decomposed in two steps. 
\begin{lem}\label{lem35}
Let $f_1$ and $f_2$ be two almost-simple functions. Denote by $f_1\oplus f_2$ the function \ $(w_1,w_2) \mapsto f_1(w_1)+f_2(w_2)$. Then $f_1\oplus f_2$ is also almost-simple.
\end{lem}
\Pv \ Consider the following decompositions for $f_1$ and $f_2$: $$\begin{cases}f_1=g_1+h_1 \\
f_2=g_2+h_2 \end{cases} $$
where $g_1$ and $g_2$ are simple functions of respective indexes $\iota_1$ and $\iota_2$, and $$\exists B_1,B_2 >0 \ , \ \begin{cases}\forall w_1, \parallel \nabla h_1(w_1) \parallel \leqslant B_1\\
\forall w_2, \parallel \nabla h_2(w_2) \parallel \leqslant B_2
 \end{cases}$$ 
Then, $g_1 \oplus g_2$ is a simple function of index $\iota_1 + \iota_2$, and \begin{align*} \parallel \nabla (h_1 \oplus h_2)(w_1,w_2) \parallel^2 & =  \parallel \nabla_{w_1}(h_1 \oplus h_2)(w_1,w_2) \parallel^2  \\ 
& \qquad \qquad \qquad + \parallel \nabla_{w_2}(h_1 \oplus h_2)(w_1,w_2) \parallel^2 \\
& =   \parallel \nabla_{w_1}h_1(w_1) \parallel^2 + \parallel \nabla_{w_2}h_2(w_2) \parallel^2\\ 
& \qquad \qquad \qquad \qquad \qquad \qquad \qquad \leqslant B_1^2 +B_2^2 \ .
 \end{align*} 
 In conclusion, $f_1 \oplus f_2=g_1\oplus g_2 + h_1 \oplus h_2$ is also a almost-simple function. \qed
 
 \begin{lem}\label{lem35bis}
Let $f_1$ and $f_2$ be two almost-simple functions. Then 
$$s(f_1\oplus f_2)=s(f_1)+s(f_2) .$$
\end{lem} 

\noindent \textit{Idea of the proof.} Thanks to the previous Lemma it makes sense to speak of the selector of the sum $f_1 \oplus f_2$. To obtain the formula $ s(f_1\oplus f_2)=s(f_1)+s(f_2)$, one can use the Morse--Barannikov complexes of $f_1$ and $f_2$ -- \cite{7}, \cite{3} -- to track the min-max of $f_1 \oplus f_2$. Such a proof is proposed in \cite{thWQ}. 
\qed

 Note that the Morse--Barannikov reduction requires $f_1$, $f_2$ and $f_1 \oplus f_2$ to be excellent Morse functions -- i.e. with all critical points having different critical values. To conclude, one may use the property of continuity of the selector together with the fact that excellent Morse functions are generic.
\qed

\vspace{0.5cm}
\subsubsection{Selector and transformation $\mathbf{T}$} \ Unfortunately, the notion of gf AS is not stable by transformation $\T$. The class of generating functions must be refined.\\

Let $g$ be a simple function of Morse index $\iota$ defined on $\R^n$. A generating function for the transformation $\T$ of $j^1 g$ is $$g_{\T}: \ (q,v) \mapsto qv-g(v)$$
It satisfies: \begin{align*} 
&1. \ \forall q \in \R^n, \nabla_v g_{\T}(q,.) \text{ is a diffeomorphism of } \R^n\\
& 2. \ \nabla g_{\T} \text{ is a diffeomorphism of } \R^n \times \R^n.
\end{align*}
In other words, $g_{\T}(q,.)$ is a simple function for all $q$, and $g_{\T}$ is \textit{globally} a simple function.

Recall that $g_{\T}$ realizes $\T(j^1 g)$ which is the $1$-graph of the Legendre transform of $g$,
$$g^{\mathbf{t}}: q \rightarrow q(\nabla g)^{-1}(q) - g((\nabla g)^{-1}(q)) . $$ 
It is also a simple function of Morse index $\iota$.

\begin{df}\label{def311}
A generating function $G$ is \textbf{globally simple} (gf \textbf{GS}) if:  \begin{align*}
&1. \ \forall q \in \R^n, \nabla_w G(q,.) \text{ is a diffeomorphism of } \R^k\\
&2. \ \nabla G \text{ is a diffeomorphism of } \R^n \times \R^k.\\
\end{align*}
\end{df}

\begin{lem}\label{lem36}${}$\\
1. If $G$ is a gf GS, then $G_{\T}$, is also a gf GS.\\
2. If $G$ is a gf GS, then it realises the $1$-graph of a simple function $g$, and $G_{\T}$ realises the $1$-graph of $g^{\mathbf{t}}$.
\end{lem}

\Pv \ 1. Let us compute the partial derivative vector $\nabla_{(v,w)}$ of $G_{\mathbf{T}}(q,v,w)$: $$\nabla_{(v,w)} G_{\T}(q;v;w)=\begin{pmatrix}
\nabla_v G_{\T}(q,v,w)\\ 
\nabla_w G_{\T}(q,v,w)\end{pmatrix}=\begin{pmatrix}
q-\nabla_v G(v,w)\\ 
-\nabla_w G(v,w)\end{pmatrix}=\begin{pmatrix}
q\\
0
\end{pmatrix} - \nabla G \ .$$ 

$\nabla G$ is a diffeomorphism, so is $\nabla_{(v,w)} G_{\T}(q,.,.)$ for all $q\in \R^n$.\\

Then, let us compute the global derivative vector $\nabla$ of $G_{\T}$~: $$\nabla G_{\T}(q;v;w)=\begin{pmatrix}
\nabla_q G_{\T}(q,v,w)\\
\nabla_v G_{\T}(q,v,w)\\ 
\nabla_w G_{\T}(q,v,w)\end{pmatrix}=\begin{pmatrix}
v\\
q-\nabla_v G(v,w)\\ 
-\nabla_w G(v,w)\end{pmatrix} \ . $$ 

Let $(Q,V,W)$ such that $$\begin{pmatrix}
Q\\
V\\
W
\end{pmatrix}=\begin{pmatrix}
v\\
q-\nabla_v G(v,w)\\ 
-\nabla_w G(v,w)\end{pmatrix} \ . $$

\noindent Necessarily $v=Q$. Then, as $\nabla_w G(Q,.)$ is a diffeomorphism, it follows that $w=(\nabla_w G(Q,.))^{-1}(W)$, and finally $q=V+ \nabla_v G(Q,w)$. This builds a inverse application which is smooth, so we conclude that $\nabla G_{\T}$ is a diffeomorphism. \\

\noindent 2. For all $q \in \R^n$, we denote by $w(q)$ the unique critical point of $G(q,.)$. Because of the hypothesis 1. in Definition \ref{def311}, the application $\nabla_w G$ is submersive onto $0$. Thus the set $\lbrace (q,w(q)) , q\in \R^n\rbrace$ is smooth, and so is the application $q \mapsto w(q)$. \\
$G$ generates the $1$-graph of the function $g$~: \ $q \mapsto G(q,w(q)) \ .$\\

\noindent Since $\nabla_w G (q,w(q))=0$, we have $\forall q \in \R^n \ , \  \nabla g (q)=\nabla_q G  (q,w(q)) \ .$\\

\noindent Let $Q$ be in $\R^n$. We must write $$Q =\nabla g (q) \  \Leftrightarrow  \ (q,w(q))=(\nabla G)^{-1}(Q,0) \ ,$$ and find the inverse application of $\nabla g$ as the projection on $\R^n$ of the inverse map $\nabla G^{-1}$ restricted to $\R^n \times \lbrace 0 \rbrace$.  As it is smooth, it proves that $\nabla g$ is a diffeomorphism.\\
This together with the first point of the Lemma ensures that $G_{\T}$ also realises the $1$-graph of a function. For any $Q$ in $\R^n$, let us write \begin{align*}
(\nabla G_{\T})^{-1}(Q,0,0)=(q,v,w) \ & \Leftrightarrow \ v=Q \ , \ q=\nabla_q G(v,w) , \ \nabla_w G(v,w) =0 \\
& \Leftrightarrow \ q=\nabla_q G(Q,w(Q)) , \ (v,w)=(Q,w(Q)) \\
& \Leftrightarrow \ q=\nabla g(Q) , \ (v,w)=(Q,w(Q))
 \end{align*}
We conclude that the $G_{\T}$ realises the $1$-graph of: $$ q \mapsto qQ- G(Q,w(Q)) \ ,  \text{ with } Q=(\nabla q)^{-1}(q) \ ,$$
which is $g_{\mathbf{t}}$. 
 
\qed \\

\begin{df}\label{def312}
A generating function $F$ defined on $\R^n \times \R^k$ is \textbf{globally almost-simple} (gf \textbf{GAS}) if it can be decomposed as a sum $$F=G+H \ ,  \text{ with } \begin{cases} G \text{ is a gf GS}\\
\exists B >0 \ , \ \parallel \nabla H(q,w) \parallel \leqslant B , \forall (q,w) \in \R^n \times \R^k \end{cases} $$
\end{df} 

\begin{rk}
A gf GAS is in particular a gf AS. 
\end{rk}

\begin{prop}
If $F$ is a gf GAS, then $F_{\T}$ is also a gf GAS. So in particular, $F_{\T}$ admits a selector.
\end{prop}

\Pv \ It is essentially the same as we have done to prove 1. in Lemma \ref{lem35} and 1. of Lemma \ref{lem36}. \qed \\

\vspace{0.5cm}
\subsubsection{Selector, sum, transformation $\mathbf{T}$ and convolution}\label{324} To work at the same time with sum and transformation $\T$ (or sum and convolution), the class of generating functions should be refined further. Indeed, the class of gf GAS used in subsection 3.2.3. for transformation $\T$ is not stable for the sum in general.

It works however if we restrict the study to the class of \textit{almost-convex} functions (Definition \ref{def37}) -- an analogue could be done with \textit{almost-concave} functions. Using Lemma \ref{lem33}, one naturally recover the definitions of the Legendre--Fenchel transform and the infimal-convolution in the following sense: 

\begin{cor}\label{cor31}${}$\\
1. Let $f$ be a almost-convex function. Then $$s(f_\mathbf{T})=f^* \ \text{ and } \ s(f_{\T \T})=f \ . $$
2. Let $f_1$ and $f_2$ be two almost-convex functions. Then
$$ s(f_{1 \s 2 })=f_1 \square f_2.$$
\end{cor}

\Pv \ \textit{1.} Let us write $f=g+h$ with $g$ simple of Morse index $\iota=0$ and $\nabla h $ bounded. Then $$f_{\T}(q,v)=qv-g(v)-h(v),$$
and for each $q$, $f_{\T}(q,.)$ is the sum of $v\mapsto qv-g(v)$ almost-concave and $-h$ which has its gradient bounded.\\
\textit{2.} If $f_1$ and $f_2$ are decomposed as $g_1+h_1$ and $g_2+h_2$, then for all $q \in \R^n$, $f_{1 \s 2}(q,.)$ is decomposed into: $$ f_{1 \s 2}(q,v)=(g_1(v)+g_2(q-v)) + (h_1(v)+h_2(q-v)) \ . $$
The fact that $\nabla h_1$ and $\nabla h_2$ are bounded implies that $$\forall q \ , \ \nabla (h_1(v)+h_2(q-v)) \text{ is bounded.}$$ \\
Moreover, $\forall q$, the Hessian matrix for $v \mapsto g_1(v)+g_2(q-v)$ at $v$ is: $$\text{Hess}_v(g_1(.)+g_2(q-.))=\text{Hess}_v(g_1)+\text{Hess}_{q-v}(g_2) \ .$$  
As the sum of two definite positive matrices is also a definite positive matrix, we obtain that $v \mapsto g_1(v)+g_2(q-v)$ is almost-simple of Morse index 0. Thus for all $q$ , $f_{1 \s 2 }(q,.)$ is also almost-convex, and the conclusion follows from Lemma \ref{lem33}. $\qed$

\begin{figure}[h]
\begin{center}
 \includegraphics[scale=0.5]{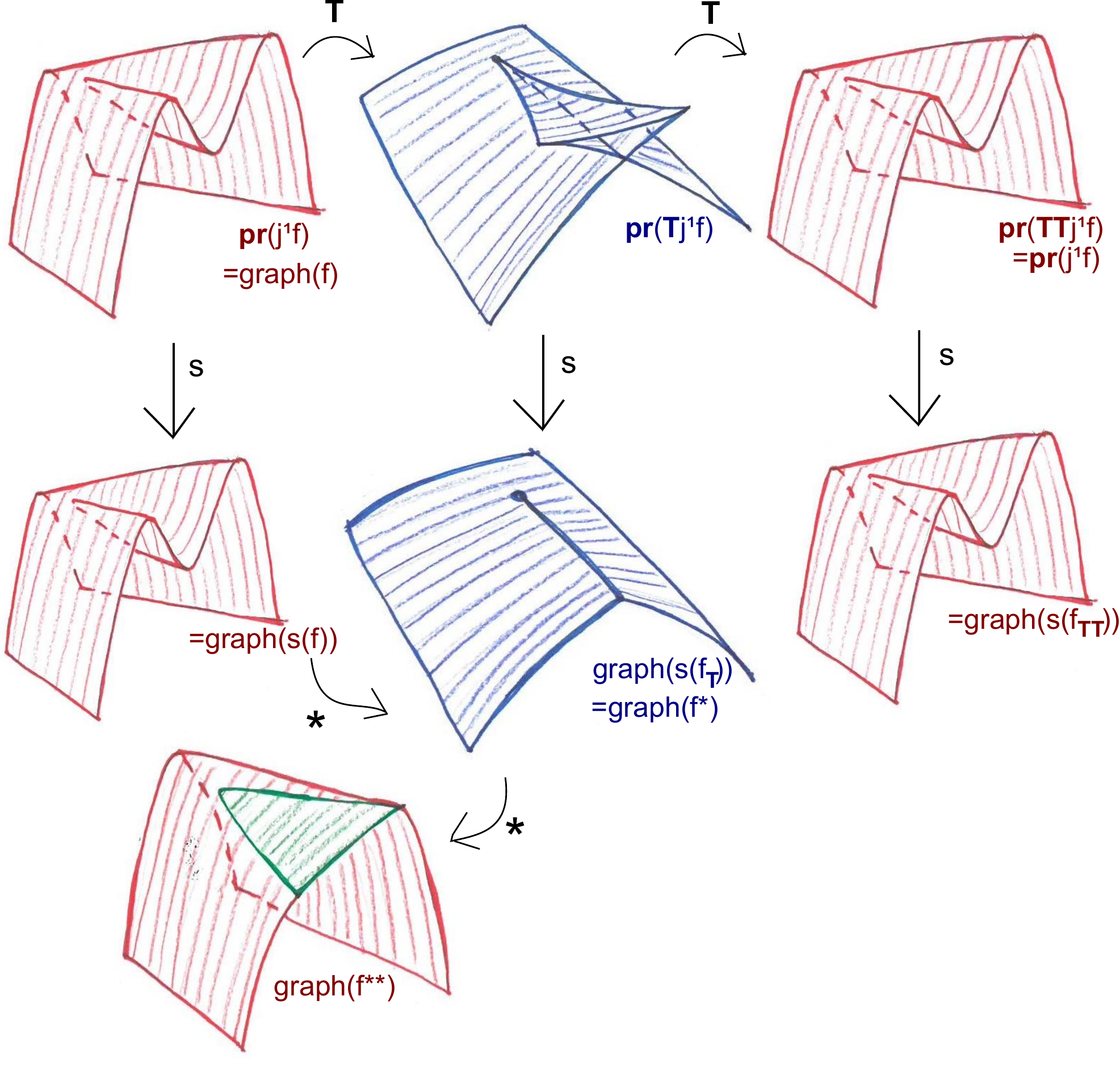} 
\end{center}
\caption{}
\label{fig9}
\end{figure}

\begin{rk}
Let $f$ be an almost-convex function, note that $s(f_{\T \T})$ does not coincide with $f^{**}$. The first provides an involution while the second convexifies $f$ (see Figures \ref{fig8} and \ref{fig9}).
\end{rk}

\subsection*{Fourier-type identities} The identities that link sum, convolution and transformation $\T$ at the geometrical level of Legendrian submanifold (identities $\underline{i})$ and $\underline{ii})$ Theorem \ref{thm21}) have an echo in convex analysis with Legendre--Fenchel transform and infimal-convolution at the level of functions: $$i) \ f_1^* + f_2^*=(f_1 \square f_2)^* \qquad \qquad ii) \ f_1^* \square f_2^*=(f_1+f_2)^*,$$ for $f_1$ and $f_2$ admissible functions -- typically almost-convex.

It is natural to suspect that $i)$ and $ii)$ correspond to the persistence of the Fourier-type identities for Legendrian submanifolds at the level of the selector of generating functions. 

\begin{thm}\label{thm327}
Let $F_1$ and $F_2$ be two gf such that $F_{(\mathbf{T}1) \p (\mathbf{T}2)}$, $F_{\mathbf{T}(1 \s 2))}$, $F_{(\mathbf{T}1) \s (\mathbf{T}2)}$ and $F_{\mathbf{T}(1 \p 2)}$ are gf AS. Then 
$$i) \ s(F_{(\mathbf{T}1) \p (\mathbf{T}2)}) = s(F_{\mathbf{T}(1 \s 2)}),$$
$$ii) \ s(F_{(\mathbf{T}1) \s (\mathbf{T}2)})=s(F_{\mathbf{T}(1 \p 2)}).$$ 
\end{thm}

\Pv \ Changing a gf by composing it at the source with a fiber-preserving diffeomorphism keeps the same $q$-parametrized Morse dynamic. In the same way, a gf is not fundamentally changed if one adds a non-degenerate quadratic form with extra variables. In terms of generating function theory, those gf are said to be \textit{equivalent} \cite{14}, \cite{4}. It is straightforward that equivalent gf's have identical selectors. 

Remind the comparison $i)$ done in Remark \ref{rk34}. It shows that $F_{(\mathbf{T}1) \p (\mathbf{T}2)}$ and $F_{\mathbf{T}(1 \s 2)}$ only differ by a fiber-preserving diffeomorphism: $\forall (q,v_1,v_2,w_1,w_2)\in \R^n\times (\R^n \times \R^n \times\R^{k_1}\times\R^{k_2}),$  $$F_{(\mathbf{T}1) \p (\mathbf{T}2)}(q,v_1,v_2,w_1,w_2)=F_{\mathbf{T}(1 \s 2)}(q,\phi(v_1,v_2,w_1,w_2)),$$
where $\phi$ is the diffeomorphism of $\R^n \times \R^n \times\R^{k_1}\times\R^{k_2}$ defined by $$\phi(v_1,v_2,w_1,w_2)=(v_1+v_2,v_1,w_1,w_2).$$ 

Thus $F_{(\mathbf{T}1) \p (\mathbf{T}2)}$ and $F_{\mathbf{T}(1 \s 2)}$ have the same selectors.\\

Similarly, the second identity can be obtained from the comparison $ii)$ of Remark \ref{rk34}, but requires a more sophisticated argument. We use a result proved in \cite{4} by Th\'{e}ret, which states that, 

if $(F_t)_{t \in [0,1]}$ is a path of gf \footnote{Th\'{e}ret's original proof is for gf which are quadratic at infinity, but it adapts well to the gf AS case.}, such that the contour of $F_t$ for all $t\in [0,1]$ is constant, then $F_0$ and $F_1$ are equivalent. 

First, we replace $F_{\mathbf{T}(1 \p 2)}: (q,v,w_1,w_2) \mapsto qv-(F_1(v,w_1)+F_2(v,w_2)) $ defined on $\R^n \times (\R^n \times \R^{k_1} \times \R^{k_2})$ by the gf defined on $\R^n \times (\R^n \times \R^n \times \R^n \times \R^{k_1} \times \R^{k_2})$ by: $$F_{\mathbf{T}(1 \p 2)}(q,v,v',V,w_1,w_2) = qv-(F_1(v,w_1)+F_2(v,w_2)) +Vv'.$$ They differ by the addition of a non-degenerate quadratic form of the extra variable $v'$ and $V$, so they are equivalent. It permits to obtain $F_0=F_{\mathbf{T}(1 \p 2)}$ and $F_1=F_{(\mathbf{T}1) \s (\mathbf{T}2)}$ defined on the same space. Then, consider the path $(F_t)_{t\in [0,1]}$ of gf defined by 
$$F_t(q,v,v',V,w_1,w_2)=qv-(F_1((1-t)v+tv',w_1)+F_2(v,w_2)) +V(v'-tv).$$

One can check that the contour of $F_t$ for all $t \in [0,1]$ is constant. Thus we conclude that $F_0$ and $F_1$ are equivalent. \qed

\newpage

\bibliography{biblio}

\end{document}